\documentclass[a4paper,11pt]{article}

\input{diagrams}

\newcommand{\op}{\ensuremath{\mbox{\hspace{1pt}{\scriptsize op}}}}

\newcommand{\cat}[1]{\ensuremath{\mbox{\bfseries {\upshape {#1}}}}}

\newcommand{\elt}[1]{\ensuremath{\mbox{\upshape elt} \hspace{1pt} {#1}}}

\newcommand{\cl}[1]{\ensuremath{\mathcal {#1}}}
\newcommand{\bb}[1]{\ensuremath{\mathbb {#1}}}

\newcommand{\lra}{\ensuremath{\longrightarrow}}
\newcommand{\map}[1]{\ensuremath{\stackrel{{#1}}{\lra}}}

\newcommand{\sm}{symmetric multicategory{}}

\newcommand{\mcat}{multicategory }

\newcommand{\bpoint}[1]{
\begin{itemize}
\item {#1}	
\end{itemize}}

\newtheorem{theorem}{Theorem}[section]
\newtheorem{proposition}[theorem]{Proposition}
\newtheorem{corollary}[theorem]{Corollary}
\newtheorem{lemma}[theorem]{Lemma}
\newtheorem{definition}[theorem]{Definition}

\newenvironment{prf}{\vspace{2ex}\begin{sloppypar}{\noindent\upshape
{\bfseries Proof. }}} {{\hspace*{\fill}
$\Box$}\end{sloppypar}\vspace{2ex}}

\newenvironment{prfof}[1]{\vspace{2ex}\begin{sloppypar}{\noindent
\upshape{\bfseries Proof of {#1}. }}} {{\hspace*{\fill}
$\Box$}\end{sloppypar}\vspace{2ex}}

\newcommand{\numroman}{\renewcommand{\labelenumi}{\roman{enumi})}}
\newcommand{\numarabic}{\renewcommand{\labelenumi}{\arabic{enumi})}}

\newcommand{\pica}{\begin{center} \input}
\newcommand{\picz}{\end{center}}
\newcommand{\length}[1]{\setlength{\unitlength}{#1}}

\newlength{\leng}
\newlength{\fontleng}
\newcommand{\sunit}{\setlength{\unitlength}{1mm}}
\newcommand{\setleng}[1]{\setlength{\leng}{#1}
    \setlength{\unitlength}{\leng}}
\newcommand{\setunit}[1]{\setlength{\unitlength}{#1}}

\newcommand{\diagc}[3]{

\put(#1,#2){
\setlength{\unitlength}{#3} %

\put(1,2){\line(1,0){5}} %
\put(1,3){\line(1,0){5}} %
\put(1,4){\line(1,0){5}} %
\put(1,1){\line(1,0){5}} %
\put(3.5,-1.5){\line(1,1){4}} %
\put(3.5,6.5){\line(1,-1){4}}}}

\newcommand{\diagf}[3]{

\put(#1,#2){
\setlength{\unitlength}{#3} %

\put(0,0){\line(1,0){4}}  %
\put(0,0){\line(2,3){2}}  %
\put(4,0){\line(-2,3){2}}}}

\newcommand{\diagj}[3]{

\put(#1,#2){
\setlength{\unitlength}{#3} %

\put(0,0){\line(1,0){4}}  %
\qbezier(0,0)(2,2)(4,0) }}

\newcommand{\diagp}[3]{

\put(#1,#2){

\setlength{\unitlength}{#3}  %
\fontsize{#3}{20}

\put(0,0){\makebox(0,0){$\equiv$}}  %
\put(0.35,0){\makebox(0,0){$\rangle$}}}}

\newcommand{\diagq}[3]{

\put(#1,#2){

\setlength{\unitlength}{#3}  %

\qbezier(0.5,0)(2.5,0.5)(3,4)   %
\put(3,4){\makebox(0,0)[bc]{{\tiny $\Delta$}}} }}

\newcommand{\diagad}[3]{

\put(#1,#2){

\setlength{\unitlength}{#3}  %

\qbezier(0.5,0)(2.5,-0.5)(3,-4)   %
\put(3,-4){\makebox(0,0)[t]{{\tiny $\nabla$}}}}}

\newcommand{\diagal}[3]{

\put(#1,#2){
\setlength{\unitlength}{#3} %

\put(0,0){\line(1,0){4}}  %
\put(0,0){\line(2,3){2}}  %
\put(4,0){\line(-2,3){2}} %
\qbezier(0,0)(2,2)(4,0) }}

\newcommand{\diagan}[3]{

\put(#1,#2){\setlength{\leng}{#3}
\setlength{\unitlength}{2\leng} %

\diagal{0}{0}{2\leng}       %
\diagf{7.3}{0}{2\leng}         %
\diagp{6}{1.5}{3\leng} }}

\sunit

\newcommand{\onetwo}[8]{

\put(#1,#2){ \setlength{\leng}{#3} %
\setlength{\unitlength}{\leng} %

\diagj{0}{0}{6\leng} %
\put(-1,0){\makebox(0,0)[r]{#4}} %
\put(26,0){\makebox(0,0)[l]{#5}} %
\put(12,7){\makebox(0,0)[b]{#6}} %
\put(12,-2){\makebox(0,0)[t]{#7}} %
\put(12,2.5){\makebox(0,0)[c]{$\Downarrow$}}  %
\put(14,2.5){\makebox(0,0)[l]{#8}} }}

\sunit

\newcommand{\diagbl}[3]{

\put(#1,#2){\setlength{\leng}{#3}
\setlength{\unitlength}{.4\leng} %

\put(0,0){\line(1,2){10}}  %
\put(10,20){\line(1,0){20}}  %
\put(30,20){\line(1,-2){10}}  %
\put(0,0){\line(1,0){40}}}}

\newcommand{\diagbn}[3]{

\put(#1,#2){\setlength{\leng}{#3}
\setlength{\unitlength}{1.2\leng} %
\fontsize{3.5\leng}{15}

\put(5,0){\vector(1,0){10}}  %
\put(10,9){\makebox(0,0)[c]{.}}  %
\put(10,3.5){\makebox(0,0)[c]{$\Downarrow$}}}} %

\sunit

\sunit

\newcommand{\onecell}[6]{

\put(#1,#2){\setlength{\leng}{#3}
\setlength{\unitlength}{1.2\leng} %

\put(4,0){\makebox(0,0)[r]{#4}}  %
\put(16.5,0){\makebox(0,0)[l]{#5}}  %
\put(10,1.5){\makebox(0,0)[b]{#6}}  %

\put(5,0){\vector(1,0){10}}}}  %

\newcommand{\threetwo}[9]{

\put(0,0){\setlength{\leng}{1.3mm}
\setlength{\unitlength}{\leng} %

\diagbl{0}{0}{\leng}  %

\put(-2,-1){\makebox(0,0)[c]{#1}}  %
\put(2,9){\makebox(0,0)[c]{#2}}  %
\put(15,9){\makebox(0,0)[c]{#3}}  %
\put(19,-1){\makebox(0,0)[c]{#4}}  %

\put(0,5){\makebox(0,0)[c]{#5}}  %
\put(8.5,11){\makebox(0,0)[c]{#6}}  %
\put(17,5){\makebox(0,0)[c]{#7}}  %
\put(8.5,-3){\makebox(0,0)[c]{#8}}  %

\put(8.5,4){\makebox(0,0)[c]{$\Downarrow$}}  %
\put(11,4){\makebox(0,0)[c]{#9}} }} %

\newcommand{\threetwob}[9]{

\put(0,0){\setlength{\leng}{1.3mm}
\setlength{\unitlength}{\leng} %

\diagbl{0}{0}{\leng}  %

\put(-2,-1){\makebox(0,0)[c]{#1}}  %
\put(2,9){\makebox(0,0)[c]{#2}}  %
\put(15,9){\makebox(0,0)[c]{#3}}  %
\put(19,-1){\makebox(0,0)[c]{#4}}  %

\put(0,5){\makebox(0,0)[c]{#5}}  %
\put(8.5,11){\makebox(0,0)[c]{#6}}  %
\put(17,5){\makebox(0,0)[c]{#7}}  %
\put(8.5,-3){\makebox(0,0)[c]{#8}}  %

\put(8.5,4){\makebox(0,0)[c]{#9}}   }} %

\newcommand{\diagbo}[3]{

\put(#1,#2){\setlength{\leng}{#3}
\setlength{\unitlength}{0.4\leng} %

\diagp{68}{12}{6\leng}

\begin{picture}(100,40)
\put(0,0){

\begin{picture}(50,30)      %

\put(0,0){\line(1,2){10}}  %
\put(10,20){\line(1,1){10}}  %
\put(20,30){\line(1,-1){10}}  %
\put(30,20){\line(1,0){15}}  %
\put(45,20){\line(1,-2){5}}  %
\put(40,0){\line(1,1){10}}  %
\put(0,0){\line(1,0){40}}  %
\put(10,20){\line(1,0){20}}  %
\put(30,20){\line(1,-2){10}}  %

\end{picture}}

\put(80,0){
\begin{picture}(50,30)      %

\put(0,0){\line(1,2){10}}  %
\put(10,20){\line(1,1){10}}  %
\put(20,30){\line(1,-1){10}}  %
\put(30,20){\line(1,0){15}}  %
\put(45,20){\line(1,-2){5}}  %
\put(40,0){\line(1,1){10}}  %
\put(0,0){\line(1,0){40}}  %

\end{picture}}

\end{picture}}}

\newcommand{\diagboo}[3]{

\put(#1,#2){\setlength{\leng}{#3}
\setlength{\unitlength}{0.4\leng} %

\diagp{68}{12}{6\leng}

\begin{picture}(100,40)
\put(0,0){

\begin{picture}(50,30)      %

\put(0,0){\line(1,2){10}}  %
\put(10,20){\line(1,1){10}}  %
\put(20,30){\line(1,-1){10}}  %
\put(30,20){\line(1,0){15}}  %
\put(45,20){\line(1,-2){5}}  %
\put(40,0){\line(1,1){10}}  %
\put(0,0){\line(1,0){40}}  %
\put(10,20){\line(1,0){20}}  %
\put(30,20){\line(1,-2){10}}  %
\put(0,0){\line(3,2){30}} 
\put(40,0){\line(1,4){5}}

\end{picture}}

\put(80,0){
\begin{picture}(50,30)      %

\put(0,0){\line(1,2){10}}  %
\put(10,20){\line(1,1){10}}  %
\put(20,30){\line(1,-1){10}}  %
\put(30,20){\line(1,0){15}}  %
\put(45,20){\line(1,-2){5}}  %
\put(40,0){\line(1,1){10}}  %
\put(0,0){\line(1,0){40}}  %

\end{picture}}

\end{picture}}}

\newcommand{\threethree}[5]{

\put(0,0){\setlength{\leng}{1mm}
\setlength{\unitlength}{0.4\leng} %

\diagp{68}{12}{6\leng} %
\put(70,20){\makebox(0,0)[b]{#5}}

\begin{picture}(100,40)
\put(0,0){

\begin{picture}(50,30)      %

\put(0,0){\line(1,2){10}}  %
\put(10,20){\line(1,1){10}}  %
\put(20,30){\line(1,-1){10}}  %
\put(30,20){\line(1,0){15}}  %
\put(45,20){\line(1,-2){5}}  %
\put(40,0){\line(1,1){10}}  %
\put(0,0){\line(1,0){40}}  %
\put(10,20){\line(1,0){20}}  %
\put(30,20){\line(1,-2){10}}  %
\put(20,10){\makebox(0,0){#3}}  %
\put(42,12){\makebox(0,0){#2}}  %
\put(21,23){\makebox(0,0){#1}}  %

\end{picture}}

\put(80,0){
\begin{picture}(50,30)      %

\put(0,0){\line(1,2){10}}  %
\put(10,20){\line(1,1){10}}  %
\put(20,30){\line(1,-1){10}}  %
\put(30,20){\line(1,0){15}}  %
\put(45,20){\line(1,-2){5}}  %
\put(40,0){\line(1,1){10}}  %
\put(0,0){\line(1,0){40}}  %
\put(25,10){\makebox(0,0){#4}}  %

\end{picture}}

\end{picture}}}

\sunit

\newcommand{\three}[4]{
\put(#1,#2){ \setlength{\leng}{#3} %
\setlength{\unitlength}{\leng} %
\diagp{0}{0}{6\leng}  %
\put(0,4){\makebox(0,0)[b]{#4}} }}

\sunit

\sunit

\sunit

\newcommand{\twotwo}[7]{

\put(0,0){ \setlength{\leng}{1mm} %
\setlength{\unitlength}{\leng} %
\diagf{0}{0}{4\leng}  %
\put(-1,-1){\makebox(0,0)[tr]{#1}}   
\put(8,14){\makebox(0,0)[b]{#2}}  %
\put(18,-1){\makebox(0,0)[tl]{#3}}  %

\put(3,6){\makebox(0,0)[br]{#4}}   
\put(14,6){\makebox(0,0)[bl]{#5}}  %
\put(8,-2){\makebox(0,0)[t]{#6}}  %

\put(8,5){\makebox(0,0)[c]{$\Downarrow$}}  %
\put(9,5){\makebox(0,0)[l]{#7}}  %
}}

\sunit

\newcommand{\twotwob}[7]{

\put(0,0){ \setlength{\leng}{1mm} %
\setlength{\unitlength}{\leng} %
\diagf{0}{0}{4\leng}  %
\put(-1,-1){\makebox(0,0)[tr]{#1}}   
\put(8,14){\makebox(0,0)[b]{#2}}  %
\put(18,-1){\makebox(0,0)[tl]{#3}}  %

\put(3,6){\makebox(0,0)[br]{#4}}   
\put(14,6){\makebox(0,0)[bl]{#5}}  %
\put(8,-2){\makebox(0,0)[t]{#6}}  %

\put(8,5){\makebox(0,0)[c]{#7}}  %

}}

\newcommand{\assleft}[5]{

\put(0,0){\setlength{\leng}{1.3mm}
\setlength{\unitlength}{\leng} %

\put(0,0){
\setlength{\unitlength}{.4\leng} %

\put(0,0){\line(1,2){10}}  %
\put(10,20){\line(1,0){20}}  %
\put(30,20){\line(1,-2){10}}  %
\put(0,0){\line(1,0){40}}   %
\put(10,20){\line(3,-2){30}}}

\put(0,5){\makebox(0,0)[c]{#1}}  %
\put(8.5,11){\makebox(0,0)[c]{#2}}  %
\put(17,5){\makebox(0,0)[c]{#3}}  %
\put(8.5,-3){\makebox(0,0)[c]{#5}}  %

\put(8.5,3.3){\makebox(0,0)[c]{#4}}   }}

\newcommand{\assright}[5]{

\put(0,0){\setlength{\leng}{1.3mm}
\setlength{\unitlength}{\leng} %

\put(0,0){
\setlength{\unitlength}{.4\leng} %

\put(0,0){\line(1,2){10}}  %
\put(10,20){\line(1,0){20}}  %
\put(30,20){\line(1,-2){10}}  %
\put(0,0){\line(1,0){40}}   %
\put(0,0){\line(3,2){30}}}

\put(0,5){\makebox(0,0)[c]{#1}}  %
\put(8.5,11){\makebox(0,0)[c]{#2}}  %
\put(17,5){\makebox(0,0)[c]{#3}}  %
\put(8.5,-3){\makebox(0,0)[c]{#5}}  %

\put(8.5,3.3){\makebox(0,0)[c]{#4}}   }}

\newcommand{\diagcz}[1]{

\setlength{\leng}{#1}
\setlength{\unitlength}{\leng} %

\put(0,0){\line(1,0){16}}}

\newcommand{\scalecq}{

\put(0,0){
\setlength{\unitlength}{0.3\leng} %

\put(20,90){\line(0,1){20}}      %
\put(40,90){\line(0,1){20}}      %
\put(80,90){\line(0,1){20}}      %

\put(20,90){\line(1,0){60}}      %
\put(20,90){\line(1,-1){30}}     %
\put(80,90){\line(-1,-1){30}}    %
\put(50,40){\line(0,1){20}}      %
\put(60,104){\makebox[0pt]{$\cdots$}} }}

\newcommand{\scalecr}[5]{

\put(0,0){
\setunit{0.3\leng} %
\scalecq
\put(20,114){\makebox[0pt]{#1}}   %
\put(40,114){\makebox[0pt]{#2}}   %
\put(80,114){\makebox[0pt]{#3}}   %
\put(50,76){\makebox(0,0){#5}}      
\put(45,38){\makebox(0,0)[t]{#4}}  }} 

\newcommand{\scalecs}[9]{

\put(0,0){
\setunit{0.3\leng} %
\put(20,0){ \scalecq
\put(20,114){\makebox[0pt]{#4}}   %
\put(40,114){\makebox[0pt][r]{#5}}   %
\put(80,114){\makebox[0pt]{#6}}   %
\put(50,76){\makebox(0,0){#9}}      
\put(50,36){\makebox(0,0)[t]{#7}}} \setunit{0.2\leng}
\put(39.8,124){\setleng{0.67\leng} \scalecr{#1}{#2}{#3}{}{#8}}}}

\usepackage{amsfonts}
\usepackage{amssymb}
\usepackage{amsmath}
\usepackage{fancyheadings}
\usepackage{latexcad}

\begin{document}

\title{Weak $n$-categories: opetopic and multitopic foundations}
\author{Eugenia Cheng\\ \\Department of Pure Mathematics, University
of Cambridge\\E-mail: e.cheng@dpmms.cam.ac.uk}
\date{October 2002}
\maketitle

\begin{abstract}

We generalise the concepts introduced by Baez and Dolan to define
opetopes constructed from symmetric operads with a category,
rather than a set, of objects.  We describe the category of
1-level generalised multicategories, a special case of the concept
introduced by Hermida, Makkai and Power, and exhibit a full
embedding of this category in the category of symmetric operads
with a category of objects.  As an analogy to the Baez-Dolan slice
construction, we exhibit a certain multicategory of function
replacement as a slice construction in the multitopic setting, and
use it to construct multitopes.  We give an explicit description
of the relationship between opetopes and multitopes.

\end{abstract}

\setcounter{tocdepth}{3}
\tableofcontents

\section*{Introduction}
\addcontentsline{toc}{section}{Introduction}

The problem of defining a weak $n$-category has been approached in
various different ways (\cite{bd1}, \cite{hmp1}, \cite{lei1},
\cite{pen1}, \cite{bat1}, \cite{tam1}, \cite{str2}, \cite{may1},
\cite{lei7}), but so far the relationship between
these approaches has not been fully understood.  The subject of
the present paper is the relationship between the approaches given
in \cite{bd1} and \cite{hmp1}.

In \cite{bd1}, John Baez and James Dolan give a definition
of weak $n$-categories based on opetopes and opetopic
sets.  In \cite{hmp1}, Claudio Hermida, Michael Makkai and
John Power begin a related definition, based on multitopes
and multitopic sets.  In each case the definition has two
components.  First, the language for describing $k$-cells
is set up.  Then, a concept of universality is introduced,
to deal with composition and coherence.  Any comparison of
the two approaches must therefore begin at the construction
of $k$-cells, and in this paper we restrict our attention
to this process.  This, in the terminology of Baez and
Dolan, is the theory of opetopes.

In \cite{bd1}, the underlying shapes of $k$-cells are shapes called `opetopes'
by Baez and Dolan.  The starting point is the theory of (symmetric)
operads. A `slicing' process on operads is defined, which is the
means of `climbing up' through dimensions; it is eventually used
to construct $(k+1)$-cells from $k$-cells. Opetopes are
constructed from the slicing process iterated, and presheaves on
the category of opetopes are called opetopic sets. A weak
$n$-category is defined as an opetopic set with certain
properties.

In \cite{hmp1}, an analogous process is presented, with shapes called
`multitopes'.  The construction is based on multicategories in a generalised
form defined in the paper. Instead of a slicing process, the construction of a
`multicategory of function replacement' is given.  This is a more general
concept, and multitopic sets are defined directly from the iteration of this
process.  Multitopes are then defined to arise from the terminal multitopic
set, and multitopic sets are shown to arise as presheaves on the category of
multitopes.

Although the multitopic approach was developed explicitly as an analogy to the
opetopic approach, the exact relationship between the notions has not
previously been clear.  The conspicuous difference between the two approaches
is the presence in the opetopic version, and absence in the multitopic, of
symmetric actions. In this paper we make explicit the relationship between
opetopes and multitopes, showing that they are `the same up to isomorphism'.  

In fact, we do not use the definition of opetopes precisely as given in
\cite{bd1}, but rather, we develop a generalisation of the notion along lines
which Baez and Dolan began but chose to abandon, for reasons unknown to the
present author.  Baez and Dolan work with operads having an arbitrary {\em set}
of types (objects), but at the beginning of the paper they use operads having
an arbitrary {\em category} of types, before restricting to the case where the
category of types is small and discrete.  However, the construction gives many
copies of each opetope, and we need to regard these as isomorphic.  So we need a
category of objects in order to preserve this vital information.  Without it, the isomorphisms are lost and such
objects are considered to be different and in this manner the relationship
between the two approaches is destroyed.  We discuss this in more detail in
Section~\ref{over}

Thus motivated, we study the approach presented by Baez and Dolan,
but using operads with a category of objects; we refer to these as
symmetric multicategories (with a category of objects), in
accordance with the terminology of \cite{hmp1} and \cite{lei1}.

The approach presented by Hermida, Makkai and Power uses
generalised 2-level multicategories, which have `upper level' and
`lower level' objects.  As far as the construction of multitopes,
however, we have found only 1-level versions to be involved, so we
consider only these, which we refer to simply as generalised
multicategories.

The constructions of multitopes and opetopes are explicitly analogous,
so we compare them step by step as follows.

We begin, in Section~\ref{over} with an informal overview of the whole theory.
We include for completeness the theory proposed by Leinster although the formal
treatment is given in a further work.  

In Section \ref{tacmcats} we define the categories {\bfseries
SymMulticat} and {\bfseries GenMulticat}, of symmetric and
generalised multicategories respectively.  These are the underlying
theories of the two approaches.  

In Section \ref{tacxi} we construct a functor
    \[\xi:\mbox{{\bfseries GenMulticat}} \longrightarrow
\mbox{{\bfseries SymMulticat}}\]
    and show that it is full and faithful.  Given a generalised
multicategory $M$, $\xi$ acts by leaving the objects unchanged,
but adding a symmetric action freely on the arrows.  (By `free'
here we mean that the orbit of an arrow with $n$ source elements
is the size of the whole permutation group ${\mathbf S}_n$.)

Clearly not all symmetric multicategories are in the image of
$\xi$.  To be in the image, a symmetric multicategory certainly
must have a discrete category of objects (we call this {\em
object-discrete}) and a free symmetric action (we call this {\em
freely symmetric}). We show that these conditions are in fact
necessary and sufficient.  Eventually we will see that every symmetric
multicategory used in the construction is equivalent to one with these
properties.  

In Section \ref{tacslice} we examine the construction of
opetopes.  We first define and compare the slicing
processes.  Our method is as follows.  Given a morphism of
symmetric multicategories \[\phi:Q \longrightarrow \xi(M)\]
we construct a morphism \[\phi^+:Q^+ \longrightarrow
\xi(M_+)\] from the action of $\phi$.  We show that if
$\phi$ is an equivalence, then $\phi^+$ is also an
equivalence.  In particular we deduce that the functor
$\xi$ and the slicing process `commute' up to equivalence,
that is, for any generalised multicategory $M$ \[\xi(M)^+
\simeq \xi(M_+).\]

In Section \ref{tacopetopes} we apply the above constructions to
opetopes and multitopes.  Writing $I$ for the symmetric
multicategory with one object and one arrow, a $k$-dimensional
opetope is defined to be an object of $I^{k+}$, the $k$th iterated
slice of $I$.  Similarly, writing $J$ for the generalised
multicategory with one object and one arrow, a $k$-dimensional
multitope is defined to be an object of $J_{k+}$, the $k$th
iterated slice of $J$.  By the above constructions, we have for
each $k$
    \[\xi(J_{k+}) \simeq I^{k+}\]
giving a correspondence between opetopes and multitopes.

Hermida, Makkai and Power suggest that where their concept is
``concrete and geometric'' the Baez-Dolan concept is ``abstract
and conceptual''. In uniting the two approaches the reward is a
concept which enjoys the elegance of  being abstract and
conceptual while at the same time providing a concrete, geometric
description of the objects in question.

\subsubsection*{Terminology}

\renewcommand{\labelenumi}{\roman{enumi})}
\begin{enumerate}

\item Since we are concerned chiefly with {\em weak}
$n$-categories, we follow Baez and Dolan (\cite{bd1}) and omit the
word `weak' unless emphasis is required; we refer to {\em strict}
$n$-categories as `strict $n$-categories'.

\item We use the term `weak $n$-functor' for an $n$-functor where
functoriality holds up to coherent isomorphisms, and `lax' functor
when the constraints are not necessarily invertible.

\item In \cite{bd1} Baez and Dolan use the terms `operad' and
`types' where we use `multicategory' and `objects'; the latter
terminology is more consistent with Leinster's use of `operad' to
describe a multicategory whose `objects-object' is 1.

\item In \cite{hmp1} Hermida, Makkai and Power use the term
`multitope' for the objects constructed in analogy with the
`opetopes' of \cite{bd1}.  This is intended to reflect the fact
that opetopes are constructed using operads but multitopes using
multicategories, a distinction that we have removed by using the
term `multicategory' in both cases.  However, we continue to use
the term `opetope' and furthermore, use it in general to refer to
the analogous objects constructed in each of the theories.

\item We regard sets as sets or discrete categories with no
notational distinction.

\end{enumerate}

\bigskip {\bfseries Acknowledgements}

This work was supported by a PhD grant from EPSRC.  I would 
like to thank Martin Hyland and Tom Leinster for their support and
guidance.

\section{Overview}
\label{over}

In this Section we give an informal overview of the opetopic foundations for
theory of $n$-categories.  This is not intended to be a rigorous treatment, but
rather, to give the reader an idea of the `spirit' of the definition, the
issues involved in making it, and the reason (as opposed to the proof) that the
different approaches in question turn out to be equivalent.  For completeness we
include here discussion of Leinster's construction (\cite{lei1}) although the
formal account is given in a further work (\cite{che8}).  

\subsection{What are opetopes?}

The defining feature of the opetopic theory of $n$-categories is, superficially,
that the underlying shapes of cells are {\em opetopes}.  Below are some examples
of opetopes at low dimensions.

$\bullet$ \hspace{1ex} 0-opetope 
\begin{picture}(20,10) \put(10,0){.} \end{picture} 

$\bullet$ \hspace{1ex} 1-opetope 
\begin{picture}(20,10) \onecell{0}{0}{1mm}{}{}{} \end{picture} 

$\bullet$ \hspace{1ex} 2-opetopes
\begin{picture}(16,20)(2,0) \diagbn{0}{0}{1mm} \end{picture}
\begin{picture}(20,15)(-2,0) \onetwo{0}{0}{0.8mm}{}{}{}{}{} \end{picture}
\begin{picture}(25,22)(-5,0) \twotwo{}{}{}{}{}{}{} \end{picture}
\begin{picture}(25,20)(-2,0) \threetwo{}{}{}{}{}{}{}{}{} \end{picture}
\begin{picture}(20,30) $\cdots$ \end{picture}

$\bullet$ \hspace{1ex} a 3-opetope 
\begin{picture}(50,30)(15,0)
\diagbo{18}{0}{1mm}  \end{picture} 

\vspace{2ex} $\bullet$ \hspace{1ex} a 4-opetope

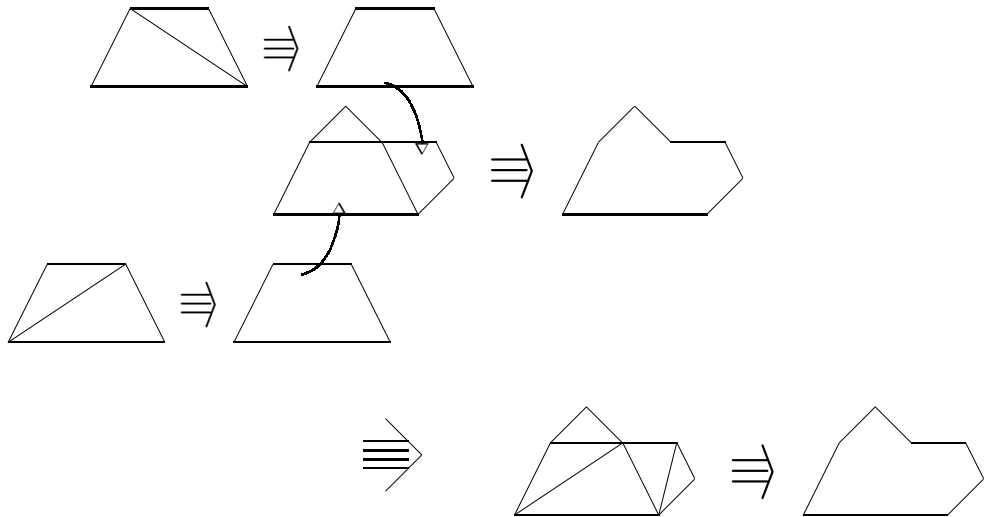
\begin{center}
\begin{picture}(100,60)
\put(26,25){\diagbo{0}{0}{1.2mm}}

\put(3,42){\assleft{}{}{}{}{}}
\put(28,47){\three{0}{0}{1mm}{}}
\put(33,42){\threetwob{}{}{}{}{}{}{}{}{}}
\put(41,42.5){\diagad{0}{0}{2mm}}

\put(-8,8){\assright{}{}{}{}{}}
\put(17,13){\three{0}{0}{1mm}{}}
\put(22,8){\threetwob{}{}{}{}{}{}{}{}{}}
\put(30,17){\diagq{0}{0}{2mm}}

\put(58,-15){\diagboo{0}{0}{1.2mm}}
\put(38,-10){\diagc{0}{0}{1.2mm}}

\end{picture} 
\end{center}
\vspace{5em}

\subsubsection*{Remarks}
\numarabic
\begin{enumerate}

\item Note that all edges and faces are directed, but we tend to omit the arrows as
at low dimensions directions are understood.

\item The number of bars on an arrow indicate its dimension.

\item The curved arrows indicate `pasting' which is otherwise difficult to
represent in higher-dimensions on a 2-dimensional sheet of paper.

\end{enumerate}

Compared with ordinary `globular' cell shapes such as
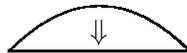
\begin{center}
\begin{picture}(35,15)(5,-3)
\onetwo{10}{0}{1mm}{$$}{$$}{$$}{$$}{$$}
\end{picture}
\end{center}
opetopes have the following important feature: the domain of a $k$-opetope is
not a single $(k-1)$-opetope but a `pasting diagram' of $(k-1)$-opetopes.  Note
that a pasting diagram can be degenerate, giving `nullary' opetopes whose
domain consists of an `empty' pasting diagram.  For example, the following is a
nullary 2-opetope:

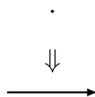
\begin{center}
\begin{picture}(20,15) \diagbn{0}{0}{1mm} \end{picture}
\end{center}

Cells in an opetopic $n$-category may thus be thought of as `labelled opetopes'
where the sources and targets of the constituent cells must match up where they
coincide in the opetope.  For example the following is a 2-cell

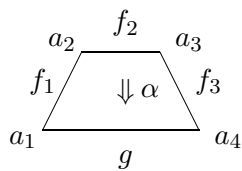
\begin{center}
\begin{picture}(25,20)(-2,-4)
\threetwo{$a_1$}{$a_2$}{$a_3$}{$a_4$}{$f_1$}{$f_2$}{$f_3$}{$g$}{$\alpha$}
\end{picture}
\end{center}
and the following is a 3-cell in which some of the lower-dimensional labels have
been omitted

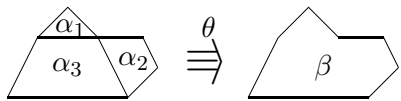
\begin{center}
\begin{picture}(35,20)
\threethree{$\alpha_1$}{$\alpha_2$}{$\alpha_3$}{$\beta$}{$\theta$}
\end{picture}
\end{center}

There now arise a philosophical question and a technical question, namely: why
and how do we do this?

\subsection{Why opetopes?}

Opetopes arise from the need, in a {\em weak} $n$-category, to record the
precise way in which a composition has been performed.  For example, consider
the following chain of composable 1-cells:
	\[a \map{f} b \map{g} c \map{h} d.\]
This gives a unique composite in an ordinary category (or any {\em strict}
$n$-category).  However, in a bicategory (or any {\em weak} $n$-category) we
should be wary of drawing such a diagram at all, as there is more than one
composite that could be produced, for example $(hg)f$ or $h(gf)$, which may in
general be distinct.

We might record the way in which the composition has occurred by a diagram such
as

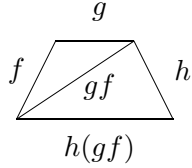
\begin{center}
\begin{picture}(25,20)(-2,-4)
\assright{$f$}{$g$}{$h$}{$gf$}{$h(gf)$}
\end{picture}
\end{center}
indicating that first $f$ is composed with $g$, and then the result is composed
with $h$.  So this diagram represents the forming of the composite $h(gf)$.  

Here the 2-cells are seen to indicate composition of their domain 1-cells.  This
is one of the fundamental ideas of the opetopic theory, that {\em composition is
not given by an operation, but by certain higher-dimensional cells}.  The cells
giving composites are those with a certain universal property, and there may be
many such cells for any composable configuration of cells.  For, as we have seen
above, there may be many distinct ways of composing a given diagram of cells.

This is the motivation behind taking opetopes as the underlying shapes of cells.

\subsection{How are opetopes constructed formally?}

We have seen that the source of a cell is to be a pasting diagram of cells
rather than just a single cell.  This is expressed using the language of {\em
multicategories}.  A multicategory is like a category whose morphisms have as
their domain a list of objects rather than just a single object.  Thus arrows
may be drawn as

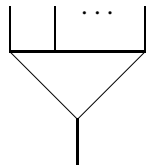
\begin{center}
\setleng{1mm}
\begin{picture}(18,30)(8,8)
\scalecq
\end{picture}
\end{center}
and composition then looks like

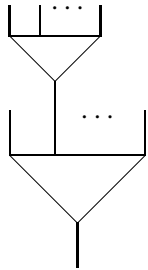
\begin{center}
\setleng{1mm}
\begin{picture}(18,55)(14,3)
\scalecs{}{}{}{}{}{}{}{}{}
\end{picture}
\end{center}

So a $k$-cell is considered as a morphism from its constituent $(k-1)$-cells to
its codomain $(k-1)$-cell.  For example

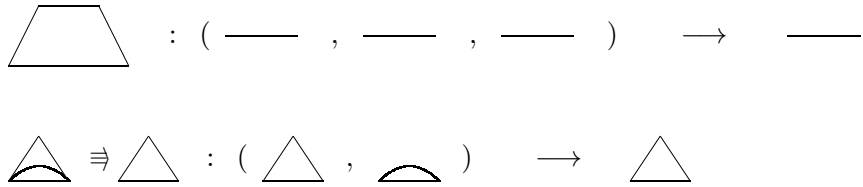
\begin{picture}(20,15)(0,3) \diagbl{0}{0}{1mm} \end{picture} : \ 
(\begin{picture}(16,10)(-2,-1) \diagcz{0.6mm} \end{picture},
\begin{picture}(16,10)(-2,-1) \diagcz{0.6mm} \end{picture},
\begin{picture}(16,10)(-2,-1) \diagcz{0.6mm} \end{picture}) \hspace{1.5em} \lra
\begin{picture}(20,10)(-8,-1) \diagcz{0.6mm} \end{picture}

\begin{picture}(25,15)(0,2) \diagan{0}{0}{1mm} \end{picture} : \ 
(\begin{picture}(13,15)(-2,2) \diagf{0}{0}{2mm} \end{picture},
\begin{picture}(13,15)(-2,2) \diagj{0}{0}{2mm} \end{picture}) \hspace{1.5em}
$\lra$ \hspace{1em} \begin{picture}(20,15)(0,2) \diagf{0}{0}{2mm}
\end{picture}
\vspace{2em}

This raises the immediate question: in what order should we list the constituent
cells?  Tom Leinster points out (\cite{lei1}) that there is no way of ordering
the cells that is stable under composition as required for
a multicategory as above.  

The three different approaches to this construction (\cite{bd1}, \cite{hmp1},
\cite{lei1}) arise from three different ways of dealing with this problem.

\bpoint{Baez and Dolan}

Baez and Dolan say: include all possible orderings.  For example

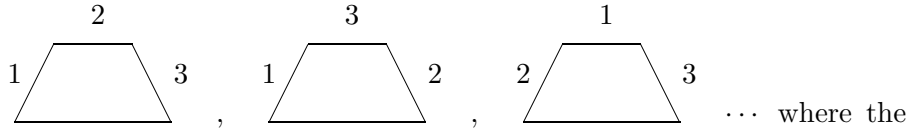
\begin{picture}(25,20)(-2,0)
\threetwob{$$}{$$}{$$}{$$}{$1$}{$2$}{$3$}{$$}{$$}
\end{picture} \ , \ 
\begin{picture}(25,20)(-2,0)
\threetwob{$$}{$$}{$$}{$$}{$1$}{$3$}{$2$}{$$}{$$}
\end{picture} \ , \ 
\begin{picture}(25,20)(-2,0)
\threetwob{$$}{$$}{$$}{$$}{$2$}{$1$}{$3$}{$$}{$$}
\end{picture}  \ $\cdots$
\vspace{3ex} where the numbers indicate the order in which the source cells are
listed.  

So a symmetric action arises, giving the different orderings, and
Baez and Dolan use {\em symmetric} multicategories for the construction.

However, a peculiar situation arise in which arrows such as

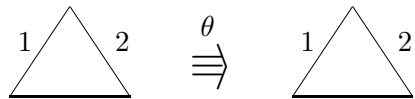
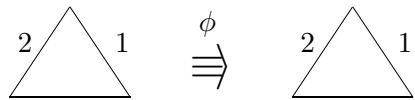
\begin{center}
\begin{picture}(25,22)(-5,-4)
\twotwob{$$}{$$}{$$}{$1$}{$2$}{$$}{$$}
\end{picture}
\begin{picture}(10,20)
\three{5}{8}{1mm}{$\theta$}
\end{picture}
\begin{picture}(25,22)(-5,-4)
\twotwob{$$}{$$}{$$}{$1$}{$2$}{$$}{$$}
\end{picture}
\end{center}
and
\begin{center}
\begin{picture}(25,22)(-5,-4)
\twotwob{$$}{$$}{$$}{$2$}{$1$}{$$}{$$}
\end{picture}
\begin{picture}(10,20)
\three{5}{8}{1mm}{$\phi$}
\end{picture}
\begin{picture}(25,22)(-5,-4)
\twotwob{$$}{$$}{$$}{$2$}{$1$}{$$}{$$}
\end{picture}
\end{center}
cannot be composed, as the ordering on the target of one does not match the
ordering of the source of the other.  The situation quickly escalates with more and more different possible
manifestations of the same opetope arising from not only the orderings on the
source cells, but also the orderings on {\em their} source cells, and so on.  
For example the following innocuous looking opetope

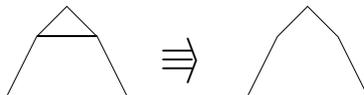
\begin{center} \length{0.4mm}

\begin{picture}(100,40)
\diagp{60}{12}{6mm}
\put(0,0){

\begin{picture}(50,30)      %

\put(0,0){\line(1,2){10}}  %
\put(10,20){\line(1,1){10}}  %
\put(20,30){\line(1,-1){10}}  %
\put(0,0){\line(1,0){40}}  %
\put(10,20){\line(1,0){20}}  %
\put(30,20){\line(1,-2){10}}  %

\end{picture}}

\put(80,0){
\begin{picture}(50,30)      %

\put(0,0){\line(1,2){10}}  %
\put(10,20){\line(1,1){10}}  %
\put(20,30){\line(1,-1){10}}  %
\put(0,0){\line(1,0){40}}  %
\put(30,20){\line(1,-2){10}}  %

\end{picture}}

\end{picture}

\end{center}
has 576 possible manifestations, and the following one

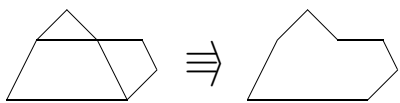
\begin{center}
\begin{picture}(40,20) \diagbo{0}{0}{1mm} \end{picture}
\end{center}
has 311040.

We need a way of saying that these objects `look the same' and this is where the
use of a {\em category} of objects comes in.  The isomorphisms in this category
tell us precisely this.

\bpoint{Hermida, Makkai, Power}

Hermida, Makkai and Power say: pick one ordering.  We know that this cannot
stable under composition; instead, the notion of multicategory is
generalised so that this stability is not required.  Rather, for each composite
there is a specified re-ordering of the source elements, satisfying some
coherence laws.  This is a notion we refer to as {\em generalised}
multicategory.

\bpoint{Leinster}

Leinster says: pick no ordering at all.  The idea is that, fundamentally,
squashing the constituent cells into a straight line is an unnatural (and indeed
rather violent) thing to try and do.  So instead, the source of an arrow such as

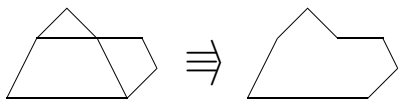
\begin{center}
\begin{picture}(40,20) \diagbo{0}{0}{1mm} \end{picture}
\end{center}
is literally the diagram

 \begin{center} \length{0.4mm}
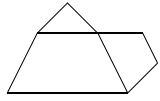
\begin{picture}(50,30)      %

\put(0,0){\line(1,2){10}}  %
\put(10,20){\line(1,1){10}}  %
\put(20,30){\line(1,-1){10}}  %
\put(30,20){\line(1,0){15}}  %
\put(45,20){\line(1,-2){5}}  %
\put(40,0){\line(1,1){10}}  %
\put(0,0){\line(1,0){40}}  %
\put(10,20){\line(1,0){20}}  %
\put(30,20){\line(1,-2){10}}  %

\end{picture}
\end{center}
expressed as a structure given by a cartesian monad $T$.  This is the notion of
$T$-multicategory. 

These differences notwithstanding, the constructions proceed in a similar
manner: a process of `slicing' is used to construct $k$-cells from $(k-1)$-cells
in each of the respective frameworks.

\subsection{Why are the different approaches equivalent?}

At first sight, it might seem implausible that a construction with so much
symmetry should give anything like a construction without any symmetry.  In
fact, the symmetric actions in the Baez-Dolan approach are a sort of {\em
trompe d'{\oe}il} created by our attempt to view constituents of an opetope in
a straight line when they simply are not in one.  It is not the opetope itself
that is symmetric, but only our presentations of it.



So, with the Baez-Dolan version, we end up with many isomorphic presentations
of the same opetope, given by all the different orders in which we could list
its components.  In effect, with the Hermida-Makkai-Power version we pick one
representative of each isomorphism class, and with the Leinster version, we take
the whole isomorphism class as one opetope.

In the end there is a trade-off between naturality (in the informal sense of
the word) and practicality.  Consider the following analogy.  If I tidy the
papers on my desk into a neat pile, I have forced them into a straight line
when they had natural positions as they were.  However, they are thus easier to
carry around.  Likewise, the Leinster construction may seem less brutal in
this sense, but the Hermida-Makkai-Power construction yields a framework that
is more practical for calculating with cells.

This means that if we are to write down a set of domain cells on a piece of
paper in a calculation, we can write them in some order.  The
Baez-Dolan construction mediates for us, giving us peace of mind that the order
we chose is irrelevant, as the symmetric actions are quietly working in the
shadows dealing with all the other possibilities.  

\subsection{How is the Baez-Dolan definition modified here?}

The present author began studying the relationship between opetopes and
multitopes as given, but began to encounter difficulties when examining the
process of slicing.  Essentially, slicing yields a multicategory whose objects
are the morphisms of the original multicategory, and whose morphisms are its
composition laws.  Given a multicategory $Q$, Baez and Dolan define the slice
multicategory $Q^+$ to be a multicategory whose {\em set} of objects is the
{\em set} of arrows of $Q$. The effect is that some information has been
abandoned, or at least, concealed. That is, we have discarded the symmetries
relating the arrows of $Q$ to one another.  As the slicing process is iterated,
progressively more information is abandoned in this manner, essentially a layer
of symmetry at each stage of slicing.  For the construction of opetopes, the
crucial fact is that the symmetries arise {\em precisely and exclusively} from the
different possible orderings of source elements.  So it is precisely these
symmetries which give the vital information about which opetopes are merely
different presentations of the same thing, and therefore should be isomorphic. 
Without it, the isomorphisms are lost and such objects are considered to be
different.  In this manner the relationship between the two approaches would
destroyed.

However, pursuing Baez and Dolan's original approach, using
multicategories with an arbitrary {\em category} of objects, it is no
longer necessary to force the category of objects of $Q^+$ to be
discrete.  This theory yields a different slice multicategory, in which the
symmetric action in $Q$ is recorded in the morphisms of the category of objects
of $Q^+$.  

This modification can then be pursued throughout the definition of $n$-category
(see \cite{che9}, \cite{che10}).  The relationship between this definition and
the original one is not currently clear.  For low dimensions it appears that
the existence of certain universal cells may eventually iron out the
differences, but such explicit arguments are unfeasible for arbitary higher
dimensions.  Moreover, such arguments cannot be applied to the structures
underlying $n$-categories where the existence of such universals has not yet
been asserted.

So what does seem clear is that the equivalences between theories as described
above facilitates much further work in this area, for example, the study of the
categories of opetopes and opetopic sets (\cite{che9}, \cite{che10}, \cite{che11},
\cite{che12}, \cite{che13}).  Using the original definition and
therefore without the help of these equivalences, this work would not have been
possible.

\section{The theory of multicategories}
\label{newopes}

Opetopes are described using the language of multicategories.  In
each of the two theories of opetopes in question, a different
underlying theory of multicategories is used.  In this section we
examine the two underlying theories, and we construct a way of
relating these theories to one another; this relationship provides
subsequent equivalences between the definitions. We adopt a concrete
approach here; certain aspects of the definitions suggest a more
abstract approach but this will require further work beyond the scope
of this work.


\label{tacmcats}

\subsection{Symmetric multicategories} \label{tacsym}
\numarabic

In \cite{bd1} opetopes are constructed using symmetric
multicategories. In this section we define {\bfseries
SymMulticat}, the category of symmetric multicategories with a
category of objects. The definition we give here includes one
axiom which appears to have been omitted from \cite{bd1}.

We write $\cl{F}$ for the `free symmetric strict monoidal
category' monad on \cat{Cat}, and ${\mathbf S}_k$ for the group of
permutations on $k$ objects; we also write $\iota$ for the
identity permutation.

\begin{definition} A {\em symmetric multicategory} $Q$ is given by
the following data
\begin{enumerate}

\item A category $o(Q)=\bb{C}$ of objects. We refer to \bb{C} as the
{\em object-category}, the morphisms of \/ ${\mathbb C}$ as {\em
object-morphisms}, and if \/ ${\mathbb C}$ is discrete, we say
that $Q$ is {\em object-discrete}.

\item For each $p\in \cl{F}{\mathbb C}^{\mbox{\scriptsize op}} \times
{\mathbb C}$, a set $Q(p)$ of arrows. Writing
    \[p=(x_1, \ldots ,x_k;x),\]
an element $f \in Q(p)$ is considered as an arrow with source and
target given by
\begin{eqnarray*}
    s(f) &=& (x_1,\ldots ,x_k)\\
    t(f) &=& x
\end{eqnarray*}
and we say $f$ has {\em arity} $k$.  We may also write $a(Q)$ for
the set of all arrows of $Q$.

\item For each object-morphism $f:x \longrightarrow y$,
an arrow $\iota(f) \in Q(x;y)$.  In particular we write
$1_x=\iota(1_x)\in Q(x;x)$.

\item Composition: for any $f \in Q(x_1, \ldots ,x_k;x)$ and
$g_i \in Q(x_{i1}, \ldots ,x_{im_i};x_i)$ for $1 \le i \le k$, a
composite
\[f \circ (g_1, \ldots ,g_k) \in Q(x_{11}, \ldots ,x_{1m_1},
\ldots ,x_{k1}, \ldots ,x_{km_k};x)\]

\item Symmetric action: for each permutation $\sigma \in {\mathbf S}_k$, a
    map
    \[\begin{array}{rccc} \sigma : & Q(x_{1}, \ldots ,x_{k};x)
    & \lra & Q(x_{\sigma (1)}, \ldots ,x_{\sigma (k)};x) \\
    & f & \longmapsto & f \sigma
    \end{array}\]

\end{enumerate}

\noindent satisfying the following axioms:

\begin{enumerate}

\item Unit laws: for any $f \in Q(x_1, \ldots ,x_m;x)$, we have
\[1_x \circ f = f = f \circ (1_{x_1}, \ldots, 1_{x_m})\]

\item Associativity: whenever both sides are defined,
 \[ \begin{array}{c} f \circ(g_1 \circ (h_{11}, \ldots , h_{1m_1}),
    \ldots , g_k \circ (h_{k1}, \ldots , h_{km_k})) = \mbox{\hspace{5em}}\\
    \mbox{\hspace*{5em}} (f \circ (g_1, \ldots , g_k)) \circ (h_{11},
    \dots ,h_{1m_1}, \ldots ,h_{k1}, \ldots ,h_{km_k}) \end{array}\]

\item For any $f \in Q(x_1, \ldots ,x_m;x)$ and
$\sigma , \sigma' \in {\mathbf S}_k$,
    \[(f \sigma)\sigma' = f(\sigma\sigma')\]

\item For any $f \in Q(x_1, \ldots ,x_k;x)$, $g_i \in Q(x_{i1},
\ldots ,x_{im_i};x_i)$ for $1 \le i \le k$, and $\sigma \in
{\mathbf S}_k$, we have
    \[(f \sigma) \circ (g_{\sigma (1)}, \ldots
    ,g_{\sigma(k)}) =
    f \circ (g_1, \ldots , g_k) \cdot \rho(\sigma)\]
where $\rho :{\mathbf S}_k \longrightarrow {\mathbf S}_{m_1 +
\ldots + m_k}$ is the obvious homomorphism.

\item For any $f \in Q(x_1, \ldots ,x_k;x)$, $g_i \in Q(x_{i1},
\ldots ,x_{im_i};x_i)$, and $\sigma_i \in {\mathbf S}_{m_i}$ for
$1 \le i \le k$, we have
    \[f \circ (g_1 \sigma_1, \ldots, g_k \sigma_k) =
    (f \circ(g_1, \ldots, g_k))\sigma\]
where $\sigma \in {\mathbf S}_{m_1 + \dots + m_k}$ is the
permutation obtained by juxtaposing the $\sigma_i$.

\item $\iota(f\circ g) = \iota(f) \circ \iota(g)$

\end{enumerate}
\end{definition}

We may draw an arrow $f \in Q(x_1, \ldots, x_k; x)$ as

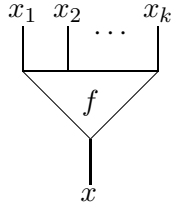
\begin{center}
\setleng{1mm}
\begin{picture}(18,30)(8,8)
\scalecr{$x_1$}{$x_2$}{$x_k$}{$\ \ x$}{$f$}
\end{picture}
\end{center} and a composite $f \circ (g_1, \ldots ,g_k)$ as

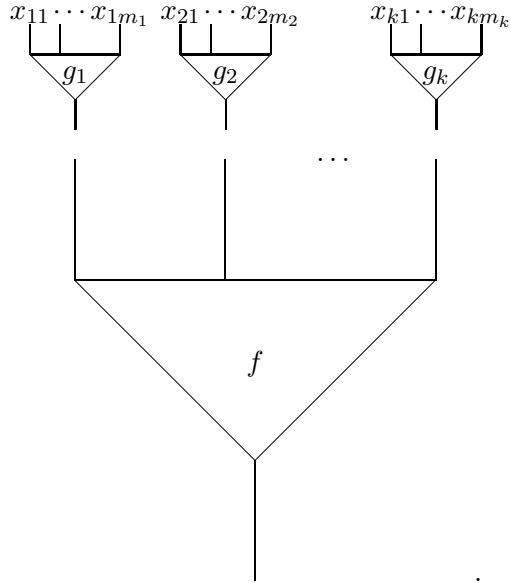
\begin{center}
\setlength{\unitlength}{0.8mm} %
\begin{picture}(70,100)(28,50)

\put(10,10){
\begin{picture}(90,140)      %
\put(20,90){\line(0,1){20}}      %
\put(45,90){\line(0,1){20}}      %
\put(63,110){\makebox(0,0){$\ldots$}}  %
\put(80,90){\line(0,1){20}}      %
\put(20,90){\line(1,0){60}}      %
\put(20,90){\line(1,-1){30}}     %
\put(80,90){\line(-1,-1){30}}    %
\put(50,40){\line(0,1){20}}      %
\put(50,76){\makebox(0,0){$f$}}      
\end{picture}}  %

\put(16,115){ \setlength{\unitlength}{0.2mm}
\begin{picture}(90,140)      %
\put(20,90){\line(0,1){20}}      %
\put(40,90){\line(0,1){20}}      %
\put(80,90){\line(0,1){20}}      %
\put(20,90){\line(1,0){60}}      %
\put(20,90){\line(1,-1){30}}     %
\put(80,90){\line(-1,-1){30}}    %
\put(50,40){\line(0,1){20}}      %
\put(20,114){\makebox[0pt]{$x_{11}$}}   %
\put(48,114){\makebox[0pt]{$\cdots$}}   %
\put(80,114){\makebox[0pt]{$x_{1m_1}$}}   %
\put(50,76){\makebox(0,0){$g_1$}}      
\put(50,36){\makebox(0,0)[t]{$$}}   
\end{picture}}

\put(41,115){ \setlength{\unitlength}{0.2mm}
\begin{picture}(90,140)      %
\put(20,90){\line(0,1){20}}      %
\put(40,90){\line(0,1){20}}      %
\put(80,90){\line(0,1){20}}      %
\put(20,90){\line(1,0){60}}      %
\put(20,90){\line(1,-1){30}}     %
\put(80,90){\line(-1,-1){30}}    %
\put(50,40){\line(0,1){20}}      %
\put(20,114){\makebox[0pt]{$x_{21}$}}   %
\put(48,114){\makebox[0pt]{$\cdots$}}   %
\put(80,114){\makebox[0pt]{$x_{2m_2}$}}   %
\put(50,76){\makebox(0,0){$g_2$}}      
\put(50,36){\makebox(0,0)[t]{$$}}   
\end{picture}}

\put(76,115){ \setlength{\unitlength}{0.2mm}
\begin{picture}(90,140)      %
\put(20,90){\line(0,1){20}}      %
\put(40,90){\line(0,1){20}}      %
\put(80,90){\line(0,1){20}}      %
\put(20,90){\line(1,0){60}}      %
\put(20,90){\line(1,-1){30}}     %
\put(80,90){\line(-1,-1){30}}    %
\put(50,40){\line(0,1){20}}      %
\put(20,114){\makebox[0pt]{$x_{k1}$}}   %
\put(48,114){\makebox[0pt]{$\cdots$}}   %
\put(80,114){\makebox[0pt]{$x_{km_k}$}}   %
\put(50,76){\makebox(0,0){$g_k$}}      
\put(50,36){\makebox(0,0)[t]{$$}}   
\end{picture}}

\end{picture}.
\end{center}

\label{tacspan} A symmetric multicategory $Q$ may be thought of as
a functor
    \[Q:\cl{F}\bb{C}^{\op} \times {\mathbb C}
    \longrightarrow \cat{Set}\]
with some extra structure.

In a more abstract view, we would expect \cl{F} to be a 2-monad on
the 2-category \cat{Cat}, which lifts via a generalised form of
distributivity to a bimonad on \cat{Prof}, the bicategory of
profunctors.  Then the Kleisli bicategory for this bimonad should
have as objects small categories, and its 1-cells should be
essentially profunctors of the form $\cl{F}\bb{C}
\makebox[0pt][l]{\hspace{1em}$\mid$}\longrightarrow \bb{D}$ in the
opposite category. However, the calculations involved in this
description are intricate and require further work.

In this abstract view, a symmetric multicategory $Q$ would then be
a monad in this bicategory.  Arrows and symmetric action (Data 2,
5) are given by the action of $Q$, identities (Data 3) by the unit
of the monad and composition (Data 4) by the multiplication for
the monad.


\begin{definition} Let $Q$ and $R$ be symmetric multicategories
with object-categories \bb{C} and \bb{D} respectively. A {\em
morphism of symmetric multicategories} $F:Q \longrightarrow R$ is
given by

\begin{itemize}

\item A functor $F=F_0:{\mathbb C} \longrightarrow {\mathbb D}$

\item For each arrow $f \in Q(x_1, \ldots ,x_k;x)$ an
arrow $Ff \in R(Fx_1, \ldots ,Fx_k; Fx)$

\end{itemize}

\noindent satisfying

\begin{itemize}

\item $F$ preserves identities: $F(\iota(f)) = \iota(Ff)$ so in
particular $F(1_x) = 1_{Fx}$

\item $F$ preserves composition: whenever it is defined
\[F(f \circ (g_1,\ldots,g_k)) = (Ff \circ (Fg_1, \ldots ,
Fg_k))\]

\item $F$ preserves symmetric action: for each $f \in Q(x_1,
\ldots , x_k;x)$ and $\sigma \in {\mathbf S}_k$
\[F(f\sigma) = (Ff) \sigma \]

\end{itemize}
\end{definition}

\noindent Composition of such morphisms is defined in the obvious
way, and there is an obvious identity morphism $1_Q:Q
\longrightarrow Q$. Thus symmetric multicategories and their
morphisms form a category {\bf SymMulticat}.

\begin{definition} A morphism $F:Q \longrightarrow R$ is an {\em equivalence} if and
only if the functor $F_0:{\mathbb C} \longrightarrow {\mathbb D}$
is an equivalence, and $F$ is full and faithful. That is, given
objects $x_1, \ldots, x_m, x$ the induced function
    \[F: Q(x_1, \ldots, x_m; x) \lra R(Fx_1, \ldots, Fx_m; Fx)\]
is an isomorphism.
\end{definition}

Note that, given morphisms of symmetric multicategories
    \[Q \map{F} R \map{G} P\]
we have a result of the form `any 2 gives 3', that is, if any two
of $F, G$ and $GF$ are equivalences, then all three are
equivalences.

Furthermore, we expect that \cat{SymMulticat} may be given the
structure of a 2-category, and that the equivalences in this
2-category would be the equivalences as above.  However, we do not
pursue this matter here.

%
%
%
%
%
%

\subsection{Generalised multicategories}
\label{tacgen}

In \cite{hmp1} multitopes are constructed using `generalised
multicategories'; in fact we need only a special case of the
generalised multicategory defined in \cite{hmp1}, that is, the
`1-level' case.

\begin{definition} A {\em generalised multicategory} $M$ is given by

\begin{itemize}

\item A set $o(M)$ of objects

\item A set $a(M)$ of arrows, with source and target
functions
\[\begin{array}{rcccl}
    s &:& a(M) & \lra & o({\mathbb C})^\star \\
    t &:& a(M) & \lra & o({\mathbb C})\\
\end{array}\]

\noindent where $A^\star$ denotes the set of lists of elements of
a set $A$. If
    \[s(f)=( x_1, \ldots ,x_k)\]
we write $s(f)_p = x_p$ and $|s(f)| = \{1,\ldots,k\}$.

\item Composition: for any $f,g \in a(M)$ with $t(g) = s(f)_p$, a
composite $f \circ_p g \in a(M)$ with
\begin{eqnarray*}
    t(f \circ_p g) & = &  t(f) \\
    |s(f \circ_p g)| & \cong & (|s(f)| \setminus \{p\}) \amalg |s(g)|
\end{eqnarray*} and {\em amalgamating maps}
\[\begin{array}{rcccc}
    \psi[f,g,p] &:& |s(f)| \setminus \{p\} & \lra & |s(f\circ_p g)|\\
    \phi[f,g,p] &:& |s(g)| & \lra & |s(f \circ_p g)|.
\end{array}\] such that $\psi \amalg \phi$ gives a bijection as above.
Equivalently, writing
\begin{eqnarray*}
    s(f) &=& ( x_1, \ldots x_k),\\
    s(g) &=& ( y_1, \ldots, y_j)
\end{eqnarray*}and
    \[(z_1, \ldots , z_{k+j-1}) = ( x_1, \ldots,
    x_{p-1}, y_1, \ldots y_j, x_{p+1}, \ldots, x_{k+j-1})\]
we have a permutation $\chi = \chi[f,g,p] \in {\mathbf S}_{k+j-1}$
such that
    \[s(f \circ_p g) = ( z_{\chi(1)}, \ldots, z_{\chi(k+j-1)}).\]

\item Identities: for each $x \in o(M)$ an arrow $1_x : x \lra x \in a(M)$

\end{itemize}

\noindent satisfying the following laws

\begin{itemize}

\item Unit laws: for any $f \in a(M)$ with $s(f)_p=x$ and
$t(f)=y$, we have
    \[\begin{array}{c} 1_y \circ_1 f = f = f \circ_p 1_x \\
    \chi[1_y,f,1]= \iota = \chi[f,1_x,p]. \end{array}\]

\item Associativity: for any $f,g,h \in a(M)$ with $s(f)_p=t(g)$ and
$s(g)_q=t(h)$ we have
    \[(f \circ_p g) \circ_{\bar{q}} h = f \circ_p(g \circ_q h)\]
where $\bar{q}=\phi[f,g,p](q)$. Furthermore, the composite
amalgamation maps must also be equal; that is, the following
coherence conditions must be satisfied:
    \[\begin{array}{c} \psi[f\circ_p g, h, \bar{q}] \circ
    \psi[f,g,p] = \psi[f, h\circ_q g, p]\\
    \psi[f\circ_p g, h, \bar{q}] \circ \bar{\phi}[f,g,p] = \phi[f,
    h\circ_q g, p] \circ \psi[g,h,q]\\
    \phi[f\circ_p g, h, \bar{q}] = \phi[f, h\circ_q g, p] \circ
    \phi[g,h,q]\end{array}\]
where $\bar{\phi}$ indicates restriction to the appropriate
domain.  Note that the conditions concern the source elements of
$f$, $g$ and $h$ respectively.

\item Commutativity: for any $f,g,h \in a(M)$ with $s(f)_p=t(g)$,\
$s(f)_q = t(h)$,\ $p \neq q$ we have
    \[(f \circ_p g)\circ_{\bar{q}} h = (f \circ_q h) \circ_{\bar{p}} g\]
    where $\bar{q}=\psi[f,g,p]$ and $\bar{p} = \psi[f,h,q]$. As above, the
composite amalgamation maps must also be equal; that is, the
following coherence conditions must be satisfied:
    \[\begin{array}{c} \psi[f\circ_p g, h, \bar{q}] \circ
    \bar{\psi}[f,g,p] = \psi[f\circ_q h, g, \bar{p}] \circ
    \bar{\psi}[f,h,q]\\
    \psi[f\circ_p g, h, \bar{q}] \circ \phi[f,g,p] =
    \phi[f\circ_q h, g, \bar{p}]\\
    \phi[f\circ_p g, h, \bar{q}] = \psi[f\circ_q h, g, \bar{p}]
    \circ \phi[f,h,q]. \end{array}\]
The conditions concern the source elements of $f$, $g$ and $h$
respectively.

\end{itemize}
\end{definition}

Note that the coherence conditions are necessary in case of
repeated source elements.


\begin{definition} A {\em morphism of generalised multicategories}
\[F=(F,\theta):M \longrightarrow N\] is given by:

\begin{itemize}

\item for each object $x \in o(M)$ an object $Fx \in
o(N)$

\item for each arrow \[f:(x_1, \ldots ,x_k)\lra x \in a(M)\] a
{\em transition map} $\theta_f= \theta_f^F \in {\mathbf S}_k$ and
an arrow \[Ff:(Fx_{\theta^{-1}(1)}, \ldots, Fx_{\theta^{-1}(k)})
\lra Fx \in a(N)\]
\end{itemize}

\noindent satisfying

\begin{itemize}

\item $F$ preserves identities: $F(1_x)=1_{Fx}$

\item $F$ preserves composition: if \(f,g \in a(M)\) and \(t(g) = s(f)_p\)
 then \[Ff \circ_{\theta_f(p)}Fg = F(f \circ_p g).\]
Furthermore, the following coherence conditions must be satisfied:
    \[\begin{array}{c} \theta_{f\circ_p g} \circ \phi[f,g,p] =
    \phi[Ff, Fg, \theta_f(p)] \circ \theta_g \\
    \theta_{f\circ_p g} \circ \psi[f,g,p] =
    \psi[Ff, Fg, \theta_f(p)] \circ \bar{\theta}_f \end{array}\]
on the source elements of $g$ and $f$ respectively, where
$\bar{\theta}$ indicates the restriction of $\theta$ as
appropriate.

\end{itemize}
\end{definition}

Given morphisms of generalised multicategories $M
\stackrel{F}{\longrightarrow} N \stackrel{G}{\longrightarrow} L$
we have a composite morphism $H = G \circ F : M \longrightarrow L$
where $H$ is the usual composite on objects and arrows, and we put
$\theta_f^H =\theta_{Ff}^G \circ \theta_f^F$. There is an identity
morphism $1_M:M \longrightarrow M$ which is the usual identity on
objects and arrows, with $\theta_f = \iota$ for all $f \in a(M)$.

Thus generalised multicategories and their morphisms form a category
{\bfseries GenMulticat}.  We now compare the two theories of
multicategories.

\subsection{Relationship between symmetric and generalised
multicategories} \label{tacxi}

We compare symmetric and generalised multicategories by means of a
functor
    \[\xi: \cat{GenMulticat} \lra \cat{SymMulticat}.\]
Given a generalised multicategory $M$, the idea is to generate a symmetric
action freely by adding in symmetric copies of each morphism.  The arrows of $M$
are then representatives of symmetry classes of arrows of $\xi(M)$.      
    
We begin by constructing this functor, and then show that it is
full and faithful.

We construct the functor $\xi$ as follows. Given a generalised
multicategory $M$, we define an object-discrete symmetric
multicategory $\xi(M) = Q$ by


\begin{itemize}

\item Objects: $o(Q)={\mathbb C}$ is the discrete category with objects
$o(M)$.

\item Arrows: for each
    \[p =(x_1, \ldots ,x_k;x) \in \cl{F}({\mathbb C})^
    {\mbox{\scriptsize op}} \times {\mathbb C}\]
an element of $Q(p)$ is given by $(f,\sigma)$ where $\sigma \in
{\mathbf S}_k$ and
    \[f : (x_{\sigma(1)}, \ldots , x_{\sigma(k)})\lra x \in
    a(M).\]

\item Composition: by commutativity, it is sufficient to define
    \[\alpha \circ_p \beta = \alpha \circ (1_{x_1}, \ldots ,
    1_{x_{p-1}}, \beta, 1_{x_{p+1}}, \ldots , 1_{x_k})\]
where
\begin{eqnarray*}
    \alpha &=& (f,\sigma) \in Q(x_1, \ldots , x_k;x)\\
    \mbox{and \ \ } \beta &=& (g,\tau) \in  Q(y_1, \ldots
    ,y_j;x_p).
\end{eqnarray*} Now given such $\alpha$ and $\beta$, we have in $M$ arrows
\[\begin{array}{rcccc}
    f &:& ( x_{\sigma(1)}, \ldots , x_{\sigma(k)}) &\lra& x\\
    \mbox{and \ \ }g &:& ( y_{\tau(1)}, \ldots ,y_{\tau(j)}) &\lra& x_p
\end{array}\]giving a composite in $M$
    \[f \circ_{\bar{p}} g: ( z_{\chi(1)}, \ldots , z_{\chi(k+j-1)})\longrightarrow x\]
where $\bar{p} = \sigma^{-1}(p)$, $\chi = \chi(f,g,\bar{p})$ and
    \[(z_1 , \ldots , z_{k+j-1}) = (x_{\sigma(1)}, \ldots ,
    x_{\sigma(\bar{p}-1)}, y_{\tau(1)}, \ldots , y_{\tau(j)} ,
    x_{\sigma(\bar{p}+1)}, \ldots, x_{\sigma(k)}).\]

We seek a composite in $Q$ with source
    \[( a_1, \ldots , a_{k+j-1}) = ( x_1, \ldots , x_{p-1},
y_1, \ldots, y_j, x_{p+1}, \ldots , x_k)\] so the composite should
be of the form $(f \circ_{\bar{p}} g, \gamma)$, where $f
\circ_{\bar{p}} g$ has source
    \[(a_{\gamma(1)}, \ldots, a_{\gamma(k+j-1)})\]
in $M$.  So we define a permutation $\gamma \in {\mathbf
S}_{j+k-1}$ by $a_{\gamma(i)} = z_{\chi(i)}$ and we define the
composite to be
    \[(f,\sigma) \circ_p (g,\tau) = (f \circ_{\bar{p}} g
    ,\gamma).\]
Note that $\gamma$ is determined by $\sigma$, $\tau$ and $\chi$.

\item For each $x \in {\mathbb C} = o(M)$, $1_x \in
Q(x;x)$ is given by $(1_x, \iota)$.

\item For each permutation $\sigma \in {\mathbf S}_k$, we have a
map
    \[\begin{array}{rccc}
    \sigma : & Q(x_1, \ldots , x_k;x) & \lra & Q(x_{\sigma(1)},
    \ldots , x_{\sigma(k)};x)
    \\
    & (f,\tau) & \longmapsto & (f,\sigma^{-1}\tau) \\
    \end{array}.\]
Note that $f$ has source $(x_{\tau(1)} \ldots, x_{\tau(k)})$ in
$M$, and $(f, \sigma^{-1}\tau)$ on the right hand side exhibits
the $i$th source of $f$ to be $x_{\sigma(\sigma^{-1}\tau)(i)} =
x_{\tau(i)}$ as required.

\end{itemize}

\numarabic

We check that this definition satisfies the conditions for a
symmetric multicategory:
\begin{enumerate}

\item Unit laws follow from unit laws of {\bf GenMulticat}

\item Associativity follows from associativity in {\bf
GenMulticat} and the coherence conditions for amalgamating maps

\item $((f,\tau)\sigma)\sigma' = (f,\sigma^{-1}\tau)\sigma' = (f,
{\sigma'}^{-1}\sigma^{-1}\tau) = (f,\tau)(\sigma\sigma')$

\item  Given
\begin{eqnarray*}
    (f,\tau) &\in& Q(x_1, \ldots , x_k;x),\\
    (g,\mu) &\in&  Q(y_1, \ldots , y_j, x_p)
\end{eqnarray*}
and $\sigma \in {\mathbf S}_k$ we check that
    \[(f,\tau)\sigma\circ_{\bar{p}} (g,\mu) =
    ((f,\tau)\circ_p(g,\mu)) \cdot \rho(\sigma)\]
where $\bar{p}=\sigma^{-1}(p)$ and $\rho$ is the homomorphism
indicated in Section~\ref{tacsym}. The required result then
follows by simultaneous composition. Note that it is sufficient to
check that both expressions in question have the same first
component and source (in $Q$), so we write $\gamma,\gamma'$ for
the permutations in the second component, without specifying what
they are. Now
    \[(f,\tau)\sigma \circ_{\bar{p}} (g,\mu) =
    (f,\sigma^{-1}\tau)\circ_{\bar{p}}(g,\mu) =
    (f\circ_{\tau^{-1}(p)}g,\gamma)\]
with source
    \[(x_{\sigma(1)}, \ldots, x_{\sigma(\bar{p}-1)}, y_1,
    \ldots, y_j, x_{\sigma(\bar{p}+1)}, \ldots, x_{\sigma(k)})\]
and
    \[((f,\tau)\circ_p(g,\mu))\cdot \rho(\sigma) =
    (f\circ_{\tau^{-1}(p)}g,\gamma')\]
with source
    \[(z_{\rho\sigma(1)}, \ldots, z_{\rho\sigma(k+j-1)})\]
where
    \[( z_1, \ldots, z_{k+j-1})=( x_1, \ldots, x_{p-1}, y_1, \ldots,
    y_j, x_{p+1}, \ldots, x_k ).\]
The action of $\rho(\sigma)$ is that of $\sigma$ on the $x_i$ but
with $( y_1, \ldots, y_j )$ substituted for $x_p$. So
    \[\begin{array}{c}
    ( z_{\rho\sigma(1)}, \ldots, z_{\rho\sigma(k+j-1)} ) =
    \mbox{\hspace{45mm}}\\
    \mbox{\hspace*{20mm}}(x_{\sigma(1)}, \ldots, x_{\sigma
    (\bar{p}-1)}, y_1, \ldots, y_j, x_{\sigma(\bar{p}+1)},
    \ldots, x_{\sigma(k)} )\end{array}\] as required.

\item Given $(f,\tau)$ and $(g,\mu)$ as above, and $\sigma \in
\cat{S}_j$ we check that
    \[(f,\tau) \circ_p (g,\mu)\sigma = ((f,\tau) \circ_p
    (g,\mu))\sigma'\]
where $\sigma' \in \cat{S}_{k+j-1}$ is given by inserting $\sigma$
at the $p$th place.

Now, on the left hand side we have
    \[\begin{array}{rcl} (f,\tau) \circ_p (g,\mu)\sigma & = &
    (f, \tau) \circ_p (g, \sigma^{-1} \mu) \\
    & = & (f \circ_{\tau^{-1}(p)} g, \gamma), \end{array}\]
say, with source
    \[(x_1, \ldots, x_{p-1}, y_{\sigma(1)}, \ldots y_{\sigma(j)},
    x_{p+1}, \ldots, x_k).\]
This agrees with the right hand side.

\item Since all object-morphisms are identities, this axiom is
trivially satisfied.

\end{enumerate}

So $\xi(M)$ is a symmetric multicategory.

Next we define $\xi$ on morphisms of generalised multicategories.
Given a morphism $F:M \longrightarrow N$ in \cat{GenMulticat} we
define a morphism \[\xi F: \xi M \longrightarrow \xi N\] in
\cat{SymMulticat} as follows.

\begin{itemize}

\item On objects: given $x \in o(\xi M) = o(M)$, put \[(\xi F)(x) =
Fx \in o(N) = o(\xi N)\]

\item On arrows: given $(f,\sigma) \in \xi M(x_1, \ldots, x_k;x)$,
put
    \[\xi F(f,\sigma) = (Ff,\sigma {\theta_f}^{-1})\]
and check that
    \[(Ff,\sigma {\theta_f}^{-1}) \in \xi N (Fx_1, \ldots, Fx_k; Fx).\]
First note that
    \[t(Ff,\sigma {\theta_f}^{-1}) = t(Ff) = F(t(f))=Fx.\]
Now
    \[s(f) = ( x_{\sigma(1)}, \ldots, x_{\sigma(k)}) \]
in $M$, so by the action of $(F,\theta)$ we have
    \[s(Ff)=( Fx_{\sigma{\theta_f}^{-1}(1)}, \ldots,
    Fx_{\sigma{\theta_f}^{-1}(k)})\]
in $N$, and so
    \[(Ff,\sigma {\theta_f}^{-1}) \in \xi N (Fx_1, \ldots, Fx_k; Fx)\]
as required.

\end{itemize}

We check that this definition satisfies the laws for a morphism of
symmetric multicategories:

\begin{itemize}

\item $\xi F$ preserves identities: since $\theta_{1_x} \in {\mathbf S}_1 =
\{\iota\}$, we have \[\xi F(1_x, \iota) = (F(1_x), \iota) =
(1_{Fx}, \iota).\]

\item $\xi F$ preserves composition: we check that
$\xi F(\alpha \circ_p \beta) = \xi F \alpha \circ_p \xi F \beta$,
and the result then follows by simultaneous composition. Put
\begin{eqnarray*}
    \alpha &=& (f,\sigma) \in Q(x_1,\ldots,x_k ; x)\\
    \mbox{and \ \ } \beta &=& (g,\tau) \in Q(y_1, \ldots , y_j ;
    y).
\end{eqnarray*} Then
\begin{eqnarray*}
    \xi F(\alpha \circ_p \beta) & = & \xi F(f \circ_{\sigma^{-1}(p)}g\
    ,\
\gamma) \\
    & = & (\ F(f \circ_{\sigma^{-1}(p)}g)\ ,\ \gamma
{\theta_f}^{-1}\ ) \\
    & = & (\ Ff \circ_{\theta_f \sigma^{-1}(p)} Fg\ ,\ \gamma
{\theta_f}^{-1}\ )
\end{eqnarray*}
and this has source
    \[s(F \alpha \circ_p F \beta)=(Fx_1, \ldots , Fx_{p-1}, Fy_1, \ldots, Fy_j,
    Fx_{p+1}, \ldots , Fx_k) .\]
For the right hand side, we have
    \[\xi F \alpha = (Ff, \sigma \theta_f^{-1})\]
    \[\xi F \beta = (Fg, \tau \theta_g^{-1})\]
and so the first component of $\xi F \alpha \circ_p \xi F \beta$
is also $Ff \circ_{\theta_f \sigma^{-1}(p)} Fg$.  So since $\xi
F(\alpha \circ_p \beta)$ and $ \xi F \alpha \circ_p \xi F \beta$
agree in the first component and source, we have the result
required.

\item $\xi F$ preserves symmetric action:
\begin{eqnarray*}
    \xi F (\ (f,\tau)\sigma\ ) & = & \xi F(f,
\sigma^{-1}\tau) \\
    & = & (Ff\ ,\ \sigma^{-1}\tau{\theta_f}^{-1}) \\
    & = & (Ff\ ,\ \tau{\theta_f}^{-1})\sigma \\
    & = & (\xi F(f,\tau))\sigma
\end{eqnarray*}

\end{itemize}

\noindent So $\xi F$ is a morphism of symmetric multicategories.

We check that $\xi$ is functorial.  Clearly $\xi 1_M = 1_{\xi M}$.
Now consider morphisms of generalised multicategories
    \[M \map{F} N \map{G} L\]
so we need to show \[\xi(G \circ F) = \xi G \circ \xi F.\]

\begin{itemize}

\item On objects
\begin{eqnarray*}
    \xi(G \circ F)(x) & = & (G \circ F)(x) \\
    & = & (\xi G \circ \xi F)(x)
\end{eqnarray*}

\item On arrows
\begin{eqnarray*}
    \xi(G \circ F)(f,\sigma) & = & (\ (G \circ F)(f)\ ,\
\sigma({\theta^{GF}}_f)^{-1}\ ) \\
    & = & (\ GFf\ ,\ \sigma ({\theta^G}_{Ff} \circ {\theta^F}_f)^{-1}\ ) \\
    & = & (\ GFf\ ,\ \sigma({\theta^F}_f)^{-1}({\theta^G}_{Ff})^{-1}\ ) \\
    & = & \xi G(\ Ff\ ,\ \sigma({\theta^F}_f)^{-1}\ ) \\
    & = & (\xi G \circ \xi F)(f, \tau)\sigma
\end{eqnarray*}

\end{itemize}

So $\xi$ is a functor as required.

\begin{proposition}
    The functor $\xi:\cat{GenMulticat} \longrightarrow
    \cat{SymMulticat}$ is full and faithful.
\end{proposition}


\begin{prf} Given any morphism
    \[G: \xi M \lra \xi N\]
of symmetric multicategories, we show that there is a unique
morphism
    \[H=(H,\theta):M\lra N\]
of generalised multicategories such that
    \[\xi H = G.\]

Suppose first that such an $H$ exists.
\begin{itemize}
\item On objects: for each object $x \in o(M) = o(\xi M)$ we must
have
    \[Hx = (\xi H)x = Gx.\]

\item On arrows: given an arrow $f \in  M(x_1, \ldots , x_k; x)$, we certainly have
\[\begin{array}{c}
    (f,\iota) \in  \xi M(x_1, \ldots, x_k; x)\\
    \mbox{and \ \ } G(f,\iota) = (\bar{f}, \sigma)
    \in  \xi N(Gx_1, \ldots Gx_k; Gx),
\end{array}\] say, where $\bar{f}$ is a morphism in $N$ with
source
    \[s(\bar{f}) = ( Gx_{\sigma(1)}, \ldots, Gx_{\sigma(k)} ).\]
Now $(\xi H)(f, \iota) = (Hf, \theta_f^{-1})$ but we must have
    \[\begin{array}{rcl} (\xi H)(f,\iota) & = & G(f, \iota) \\
    & = & (\bar{f}, \sigma) \end{array}\]
so we must have $Hf=\bar{f}$ and $\theta_f = \sigma^{-1}$.

\end{itemize}

So we define $H$ as above and check that this satisfies the axioms
for a morphism of generalised multicategories.

\begin{itemize}
\item $H$ preserves identities
\end{itemize}
We have
    \[G(1_x, \iota) = (1_{Gx}, \iota)\]
so
    \[H(1_x)=1_{Gx}=1_{Hx}.\]

\begin{itemize}
\item $H$ preserves composition
\end{itemize}
We need to show
    \[Hf \circ_{\theta_f(p)} Hg = H(f\circ_p g)\]
and that the coherence conditions are satisfied.  Now, $G$
preserves the composition of $\xi M$ so
    \[G\alpha \circ_p G\beta = G(\alpha \circ_p \beta).\]
Now we have
    \[\begin{array}{rcl} G\alpha \circ_p G\beta & = & (\bar{f},
    \theta_f^{-1}) \circ_p (\bar{g}, \theta_g^{-1})\\
    & = & (\bar{f} \circ_{\theta_f(p)} \bar{g}, \gamma), \mbox{\
    say}
    \end{array}\]
and
    \[\begin{array}{rcl} G(\alpha \circ_p \beta) & = & G(f\circ_p
    g, \gamma') \\
    & = & (\overline{f \circ_p g}, \gamma''), \mbox{\ say.}\end{array}\]
So these must be equal on both components.  Comparing first
components, we have
    \[\overline{f\circ_p g} = \bar{f} \circ_{\theta_f(p)} \bar{g}\]
but by definition we have
    \[\begin{array}{rcl} \overline{f \circ_p g} & = & H(f \circ_p g)\\
    \mbox{and\ \ } \bar{f}\circ_{\theta_f(p)} \bar{g} & = & Hf
    \circ_{\theta_f(p)} Hg \end{array}\]
so
    \[Hf \circ_{\theta_f(p)} Hg = H(f\circ_p g)\]
as required.  Furthermore, equality of the second components gives
precisely the coherence condition we require, since $\gamma$ is
formed from $\theta_f, \theta_g$ and the amalgamation map
$\chi(\bar{f}, \bar{g}, \theta_f(p))$, and $\gamma''$ is formed
from $\chi(f, g, p)$ and $\theta_{f \circ_p g}$.

So $H$ is a morphism of generalised multicategories; by
construction it is unique such that $\xi H = G$, so $\xi$ is
indeed full and faithful.  \end{prf}

We now give necessary and sufficient conditions for a symmetric
multicategory to be in the image of $\xi$.

\begin{definition} We say that a symmetric multicategory $ Q$ is {\em freely
symmetric} if and only if for every arrow $\alpha \in  Q$ and
permutation $\sigma$ \[\alpha \sigma = \alpha \Rightarrow \sigma =
\iota.\] \end{definition}

\begin{proposition} \label{prop119} Let $ Q$ be a symmetric multicategory. Then $ Q \cong
\xi(M)$ for some generalised multicategory $M$ if and only if $ Q$
is object-discrete and freely symmetric.
\end{proposition}


\begin{prf} Suppose $ Q \cong \xi(M)$. Then by the definition of $\xi$,
$ Q$ is object-discrete, with object-category ${\mathbb C}\cong
o(M)$. To show that $Q$ is freely symmetric, write $p=( x_1,
\ldots, x_k;x )$, so
    \[\begin{array}{rcl} Q(p) = \{(f,\tau) & | & f \in a(M), \ \tau
    \in {\mathbf S}_k\\
    && f: x_{\tau(1)}, \ldots, x_{\tau(k)} \lra x \in M \}\end{array}\]
and consider $\alpha = (f,\tau) \in Q(p)$. Now $(f,\tau)\sigma =
(f,\sigma^{-1}\tau)$ so
\begin{eqnarray*}
    \alpha\sigma=\alpha & \Rightarrow & \sigma^{-1}\tau = \tau \\
    & \Rightarrow & \sigma = \iota
\end{eqnarray*} as required.

Conversely, suppose that $Q$ is object-discrete and freely
symmetric.  So, given an arrow $\alpha$ of arity $k$, we have
distinct arrows $\alpha\sigma$ for each $\sigma\in{\mathbf S}_k$.
We define an equivalence relation $\sim$ on $a(Q)$, by
    \[\alpha \sim \beta \iff \beta = \alpha\sigma \mbox{ for some
    permutation } \sigma\]
and we specify a representative of each equivalence class.

Now let $M$ be a generalised multicategory whose objects are those
of $Q$, and whose arrows are the chosen representatives of the
equivalence classes of $\sim$. Composition is inherited, with
amalgamation maps re-ordering the sources as necessary.  So
associativity and commutativity are inherited; the coherence
conditions for amalgamation maps are satisfied since $Q$ is freely
symmetric. Observe that for each $x \in {\mathbb C}$, the
equivalence class of $1_x$ is $\{1_x\}$, so $M$ inherits
identities.

So $M$ is a generalised multicategory, and $\xi(M) \cong Q$. Note
that a different choice of representatives would give an
equivalent generalised multicategory. \end{prf}


\begin{definition} We call a symmetric multicategory {\em tidy} if it is
freely symmetric with a category of objects equivalent to a
discrete one.  We write \cat{TidySymMulticat} for the full
subcategory of \cat{SymMulticat} whose objects are tidy symmetric
multicategories. \end{definition}

\begin{lemma} \label{tidylemma}  A \sm\ is tidy if and only if it is equivalent to
one in the image of $\xi$.  \end{lemma}

\begin{prf}  We show that $Q$ is tidy if and only if $Q \simeq R$
where $R$ is freely symmetric and object-discrete.  The result
then follows by Proposition~\ref{prop119}.

Suppose $Q$ is tidy.  We construct $R$ as follows.  Let ${\mathbb
C}$ be the category of objects of $Q$, with ${\mathbb C}$
equivalent to a discrete category $S$, say, by
    \[{\mathbb C} {{F \atop \longrightarrow} \atop {\longleftarrow
    \atop G}} S.\]
Then $R$ is given by

\begin{itemize}

\item $o(R) = S$.

\item $R(d_1,\ldots,d_n;d) = Q(Gd_1,\ldots,Gd_n;Gd)$.

\item identities, composition and symmetric action induced from
$Q$.

\end{itemize}

\noindent Then certainly $Q \simeq R$ and $R$ is freely symmetric
and object-discrete; the converse is clear.  \end{prf}

We will later see (Section~\ref{tacopetopes}) that only tidy
symmetric multicategories are needed for the construction of
opetopes.  We now include another result that will be useful in
the next section.

\begin{lemma} If $Q$ is a tidy \sm\ then $\elt{Q}$ is equivalent
to a discrete category. \end{lemma}

\begin{prf} This may be proved by direct calculation; it is also
seen in Proposition~\ref{mcatpropf}. \end{prf}

Note that we write $\elt{Q}$ for the category of elements of $Q$,
where $Q$ is here considered as a functor $Q:\cl{F}{\mathbb
C}^{\mbox{\scriptsize op}} \times {\mathbb C} \longrightarrow
\mbox{\cat{Set}}$ with certain extra structure.

So $\elt{Q}$ has as objects pairs $(p,g)$ with $p \in
\cl{F}{\mathbb C}^{\mbox{\scriptsize op}} \times {\mathbb C}$ and
$g \in Q(p)$; a morphism $\alpha:(p,g) \longrightarrow (p',g')$ is
an arrow
        $\alpha:p \longrightarrow p' \in \cl{F}{\mathbb C}
        ^{\mbox{\scriptsize op}} \times {\mathbb C}$
such that
    \[\begin{array}{rccc}
    Q(\alpha) : & Q(p) & \lra & Q(p')
    \\
    & g & \longmapsto & g' \ . \\
    \end{array}\]
For example, an arrow
    \[(\sigma,f_1,f_2,f_3,f_4;f):(x_1,x_2,x_3,x_4;x) \longrightarrow
    (y_1,y_2,y_3,y_4;y) \in \cl{F}{\mathbb
    C}^{\mbox {\scriptsize op}} \times {\mathbb C}\]
may be represented by the following diagram

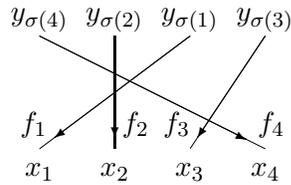
\setlength{\unitlength}{0.5mm}

\begin{center}
\begin{picture}(70,140)(20,15)
\put(0,20){ %

\begin{picture}(70,50)      %
\put(50,10){\line(0,1){30}}    %
\put(50,40){\vector(0,-1){15}}   %
\put(50,44){\makebox(0,0)[b]{$x$}} 
\put(50,6){\makebox(0,0)[t]{$y$}}  
\put(53,26){\makebox(0,0)[l]{$f$}} 
\end{picture}}

\put(0,70){

\begin{picture}(90,80)      %

\put(20,70){\line(2,-1){60}}        %
\put(20,70){\vector(2,-1){55}}      %
\put(40,70){\line(0,-1){30}}        %
\put(40,70){\vector(0,-1){27}}      %
\put(60,70){\line(-4,-3){40}}       %
\put(60,70){\vector(-4,-3){36}}     %
\put(80,70){\line(-2,-3){20}}       %
\put(80,70){\vector(-2,-3){18}}     %

\put(20,74){\makebox[0pt]{$y_{\sigma(4)}$}}   
\put(40,74){\makebox[0pt]{$y_{\sigma(2)}$}}   %
\put(60,74){\makebox[0pt]{$y_{\sigma(1)}$}}   %
\put(80,74){\makebox[0pt]{$y_{\sigma(3)}$}}   %

\put(20,36){\makebox(0,0)[t]{$x_1$}}   
\put(40,36){\makebox(0,0)[t]{$x_2$}}   %
\put(60,36){\makebox(0,0)[t]{$x_3$}}   %
\put(80,36){\makebox(0,0)[t]{$x_4$}}   %

\put(22,43){\makebox(0,0)[br]{$f_1$}}   
\put(42,43){\makebox(0,0)[bl]{$f_2$}}   %
\put(60,43){\makebox(0,0)[br]{$f_3$}}   %
\put(78,43){\makebox(0,0)[bl]{$f_4$}}   %

\end{picture}}
\end{picture}.
\end{center}Then, given any arrow $g \in Q(x_1, \ldots x_m; x)$,
we have an arrow
    \[\alpha(g) = g' \in Q(y_1, \ldots, y_m; y)\]
given by
    \[g' = (\iota(f) \circ g \circ (\iota(f_1), \ldots,
    \iota(f_m))\sigma).\]
So continuing the above example we may have:

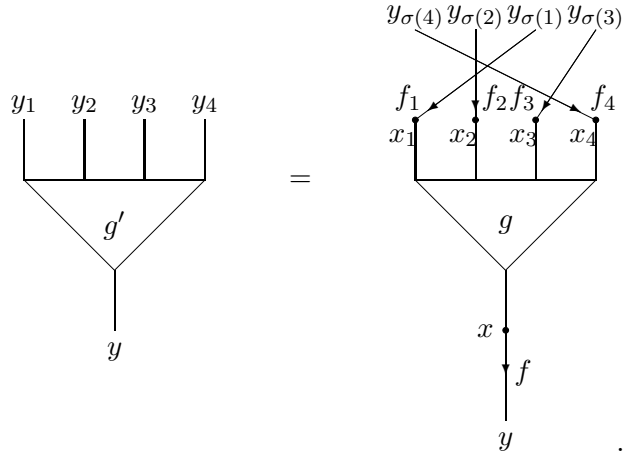
\begin{center}
\setlength{\unitlength}{0.4mm}
\begin{picture}(200,150)(40,0) 
\put(150,0){
\begin{picture}(70,50)

\put(50,10){\line(0,1){30}}    %
\put(50,40){\vector(0,-1){15}}   %
\put(46,40){\makebox(0,0)[r]{$x$}} 
\put(50,6){\makebox(0,0)[t]{$y$}}  
\put(53,26){\makebox(0,0)[l]{$f$}} 
\put(50,40){\circle*{2}}

\end{picture}}

\put(150,0){
\begin{picture}(90,140) 

\put(20,90){\line(0,1){20}}      %
\put(40,90){\line(0,1){20}}      %
\put(60,90){\line(0,1){20}}      %
\put(80,90){\line(0,1){20}}      %

\put(20,90){\line(1,0){60}}      %
\put(20,90){\line(1,-1){30}}     %
\put(80,90){\line(-1,-1){30}}    %
\put(50,40){\line(0,1){20}}      %

\put(50,75){\makebox(0,0){$g$}}      

\end{picture}}

\put(150,70){
\begin{picture}(90,80) 

\put(20,70){\line(2,-1){60}}        %
\put(20,70){\vector(2,-1){55}}      %
\put(40,70){\line(0,-1){30}}        %
\put(40,70){\vector(0,-1){27}}      %
\put(60,70){\line(-4,-3){40}}       %
\put(60,70){\vector(-4,-3){36}}     %
\put(80,70){\line(-2,-3){20}}       %
\put(80,70){\vector(-2,-3){18}}     %

\put(20,74){\makebox[0pt]{$y_{\sigma(4)}$}}   
\put(40,74){\makebox[0pt]{$y_{\sigma(2)}$}}   %
\put(60,74){\makebox[0pt]{$y_{\sigma(1)}$}}   %
\put(80,74){\makebox[0pt]{$y_{\sigma(3)}$}}   %

\put(16,37){\makebox(0,0)[t]{$x_1$}}   
\put(36,37){\makebox(0,0)[t]{$x_2$}}   %
\put(56,37){\makebox(0,0)[t]{$x_3$}}   %
\put(76,37){\makebox(0,0)[t]{$x_4$}}   %

\put(20,40){\circle*{2}} %
\put(40,40){\circle*{2}} %
\put(60,40){\circle*{2}} %
\put(80,40){\circle*{2}} %

\put(22,43){\makebox(0,0)[br]{$f_1$}}   
\put(42,43){\makebox(0,0)[bl]{$f_2$}}   %
\put(60,43){\makebox(0,0)[br]{$f_3$}}   %
\put(78,43){\makebox(0,0)[bl]{$f_4$}}   %

\end{picture}}


\put(20,0){

\begin{picture}(90,140)      %

\put(20,90){\line(0,1){20}}      %
\put(40,90){\line(0,1){20}}      %
\put(60,90){\line(0,1){20}}      %
\put(80,90){\line(0,1){20}}      %

\put(20,90){\line(1,0){60}}      %
\put(20,90){\line(1,-1){30}}     %
\put(80,90){\line(-1,-1){30}}    %
\put(50,40){\line(0,1){20}}      %

\put(20,114){\makebox[0pt]{$y_1$}}   %
\put(40,114){\makebox[0pt]{$y_2$}}   %
\put(60,114){\makebox[0pt]{$y_3$}}   %
\put(80,114){\makebox[0pt]{$y_4$}}   %

\put(50,75){\makebox(0,0){$g'$}}      
\put(50,36){\makebox(0,0)[t]{$y$}}   

\end{picture}}

\put(135,90){\makebox(0,0){$=$}}

\end{picture}.
\end{center} Note that we may write an object $(p,g)\in \mbox{elt}
(Q)$ simply as $g$, since $p$ is uniquely determined by $g$.

\section{The theory of opetopes} \label{opeope}

In this section we give the analogous constructions of opetopes in
each theory, and show in what sense they are equivalent.  That is,
we show that the respective categories of $k$-opetopes are
equivalent.

\label{tacslice}

We first discuss the process by which $(k+1)$-cells are
constructed from $k$-cells.  In \cite{bd1}, the `slice'
construction is used, giving for any \sm\ $Q$ the slice \mcat\
$Q^+$.  In \cite{hmp1} the `multicategory of function replacement'
is used but this has a more far-reaching role than that of the
Baez-Dolan slice.  For comparison with the Baez-Dolan theory, we
construct a `slice' which is analogous to the Baez-Dolan slice and
is a special case of a multicategory of function replacement.

Opetopes and multitopes are then constructed by iterating the slicing
process.  We finally apply the results already established to show
that the category of multitopes is equivalent to the category of
opetopes.  

\subsection{Slicing a symmetric multicategory}
\label{tacsymslice}

Let $Q$ be a symmetric multicategory with a category ${\mathbb C}$
of objects, so $Q$ may be considered as a functor
$Q:\cl{F}{\mathbb C}^{\mbox{\scriptsize op}} \times {\mathbb C}
\longrightarrow \mbox{\cat{Set}}$ with certain extra structure.
The slice multicategory $Q^+$ is given by:

\begin{itemize}
\item Objects: put $o(Q^+) = \mbox{elt}(Q)$
\end{itemize}

So the category $o(Q^+)$ has as objects pairs $(p,g)$ with $p \in
\cl{F}{\mathbb C}^{\mbox{\scriptsize op}} \times {\mathbb C}$ and
$g \in Q(p)$; a morphism $\alpha:(p,g) \longrightarrow (p',g')$ is
an arrow
        $\alpha:p \longrightarrow p' \in \cl{F}{\mathbb C}
        ^{\mbox{\scriptsize op}} \times {\mathbb C}$
such that
    \[\begin{array}{rccc}
    Q(\alpha) : & Q(p) & \lra & Q(p')
    \\
    & g & \longmapsto & g' \  \\
    \end{array}\]

Then, given any arrow \[g \in Q(x_1, \ldots x_m; x)\] we have an
arrow $\alpha(g) = g' \in Q(y_1, \ldots, y_m; y)$ given by
    \[g' = (\iota(f) \circ g \circ (\iota(f_1), \ldots,
    \iota(f_m))\sigma)\]
(see Section~\ref{tacxi}).

\begin{itemize}
\item Arrows: $Q^+(f_1, \ldots, f_n;f)$ is given by the set of
`configurations' for composing $f_1, \ldots, f_n$ as arrows of
$Q$, to yield $f$.
\end{itemize}

Writing $f_i \in Q(x_{i1}, \ldots x_{im_i}; x_i)$ for $1 \leq i
\leq n$, such a configuration is given by $(T,\rho, \tau)$ where

\begin{enumerate}

\item $T$ is a planar tree with $n$ nodes.  Each node is labelled
by one of the $f_i$, and each edge is labelled by an
object-morphism of $Q$ in such a way that the (unique) node
labelled by $f_i$ has precisely $m_i$ edges going in from above,
labelled by $a_{i1}, \ldots, a_{im_i} \in \mbox{arr}({\mathbb
C})$, and the edge coming out is labelled $a_i \in a({\mathbb
C})$, where $\mbox{cod}(a_{ij}) = x_{ij}$ and $\mbox{dom}(a_i) =
x_i$.

\item $\rho \in {\mathbf S}_k$ where $k$ is the number of leaves
of $T$.

\item $\tau:\{\mbox{nodes of } T\} \longrightarrow [n]=\{1, \ldots, n\}$ is a
bijection such that the node $N$ is labelled by $f_{\tau(N)}$.
(This specification is necessary to allow for the possibility $f_i
= f_j,\ i \neq j$.)

\end{enumerate} Note that $(T,\rho)$ may be considered as a `combed tree',
that is, a planar tree with a `twisting' of branches at the top
given by $\rho$.

The arrow resulting from this composition is given by composing
the $f_i$ according to their positions in $T$, with the $a_{ij}$
acting as arrows $\iota(a_{ij})$ of $Q$, and then applying $\rho$
according to the symmetric action on $Q$. This construction
uniquely determines an arrow $(T,\rho,\tau) \in Q^+(f_1, \ldots,
f_n;f)$.

\begin{itemize}
\item Composition
\end{itemize}

When it can be defined, $(T_1,\rho_1,\tau_1)
    \circ_m (T_2,\rho_2,\tau_2) = (T,\rho,\tau)$ is given by

\begin{enumerate}

\item $(T,\rho)$ is the combed tree obtained by replacing the node
${\tau_1}^{-1}(m)$ by the tree $(T_2,\rho_2)$, composing the edge
labels as morphisms of ${\mathbb C}$, and then `combing' the tree
so that all twists are at the top.

\item $\tau$ is the bijection which inserts the source of $T_2$
into that of $T_1$ at the $m$th place.

\end{enumerate}

\begin{itemize}

\item Identities: given an object-morphism
    \[\alpha=(\sigma, f_1, \ldots, f_m;f) : g \longrightarrow g',\]
$\iota(\alpha) \in Q^+(g;g')$ is given by a tree with one node,
labelled by $g$, twist $\sigma$, and edges labelled by the $f_i$
and $f$ as in the example above.

\item Symmetric action: $(T,\rho,\tau)\sigma = (T,\rho,\sigma^{-1}\tau)$

\end{itemize}

\noindent This is easily seen to satisfy the axioms for a
symmetric multicategory.

Note that, given a labelled tree $T$ with $n$ nodes and $k$
leaves, there is an arrow $(T,\rho,\tau) \in a(Q^+)$ for every
permutation $\rho \in {\mathbf S}_k$ and every bijection
$\tau:\{\mbox{nodes of } T\} \longrightarrow [n]$. Suppose
\begin{eqnarray*}
    s(T,\rho,\tau) &=& ( f_1, \ldots, f_n ) \\
    \mbox{and \ \ }t(T,\rho,\tau) &=& f.
\end{eqnarray*}
Then, given any $\rho_1 \in {\mathbf S}_k, \ \tau:\{\mbox{nodes of
} T\} \longrightarrow [n]$, we have
\begin{eqnarray*}s(T,\rho_1\rho,\tau) &=& ( f_1, \ldots, f_n ) \\
   \mbox{and \ \ }t(T,\rho_1\rho,\tau) &=& f\rho_1
\end{eqnarray*}
whereas
\begin{eqnarray*} s(T,\rho,\tau_1\tau) &=& ( f_{{\tau_1}^{-1}(1)},
    \ldots f_{{\tau_1}^{-1}(n)} ) \\
    \mbox{and \ \ }t(T,\rho,\tau_1\tau) &=& f.
\end{eqnarray*} We observe immediately that $Q^+$ is freely symmetric,
since
\begin{eqnarray*} (T,\rho,\tau)\sigma = (T,\rho,\tau) &
\Rightarrow & \sigma^{-1}\tau = \tau \\
    & \Rightarrow & \sigma = \iota.
\end{eqnarray*}
However $Q^+$ is not in general object-discrete; we will later see
(Proposition~\ref{mcatpropf}) that $Q^+$ is tidy if $Q$ is tidy.

\subsection {Slicing a generalised multicategory}
\label{tacgenslice}

Given a generalised multicategory $M$, we define a slice
multicategory $M_+$.  We use the `multicategory of function
replacement' as defined in \cite{hmp1}, which plays a role similar
to (but more far-reaching than) that of the Baez-Dolan slice.  The
slice defined in this section is only a special case of a
multicategory of function replacement, but it is sufficient for
the construction of multitopes.  Moreover, for the purpose of
comparison it is later helpful to be able to use this closer
analogy of the Baez-Dolan slice.

We first explain how this slice arises from the multicategory of
function replacement as defined in \cite{hmp1}, and then give an
explicit construction of the slice multicategory that is analogous
to the symmetric case.  This latter construction is the one we
continue to use in the rest of the work.

Using the terminology of \cite{hmp1}, the slice is defined as
follows.  Let ${\mathcal L}$ be the language with objects $o(M)$
and arrows $a(M)$, and let ${\mathbb F}$ be the free generalised
multicategory on ${\mathcal L}$.  So the objects of ${\mathbb F}$
are the objects of $M$, and the arrows of ${\mathbb F}$ are formal
composites of arrows of $M$. We define a morphism of generalised
multicategories $h:{\mathbb F} \longrightarrow M$ as the identity
on objects, and on arrows the action of composing the formal
composite to yield an arrow of M. Then we define $M_+$ to be the
multicategory of function replacement on $({\mathcal L}, {\mathbb
F}, h)$.


Explicitly, the slice multicategory $M_+$ is a generalised
multicategory given by:

\begin{itemize}

\item Objects: $o(M_+) = a(M)$.

\item Arrows: $a(M_+)$ is given by configurations for composing arrows of
$M$.

\end{itemize}

\noindent Such a configuration is given by $T=(T,\rho_T,\tau_T)$,
where: \numroman
\begin{enumerate}

\item $T$ is a planar tree with $n$ nodes labelled by $f_1, \ldots
f_n \in a(M)$, and edges labelled by objects of $M$ in such a way
that, writing
    \[s(f_i) = ( x_{i1}, \ldots, x_{im_i} ),\]
the node labelled by $f_i$ has $m$ edges coming in, labelled by $
x_{i1}, \ldots, x_{im_i}$ from left to right, and one edge going
out, labelled by $t(f_i)$.

\item $\rho_T \in {\mathbf S}_k$, where $k$ is the number of
leaves of $T$. The composition in $M$ given by $T$ has specified
amalgamation maps giving information about the ordering of the
source; $\rho_T$ is the permutation induced on the source.

\item $\tau_T:\{\mbox{nodes of }T\} \longrightarrow [n]$ is a
bijection so that the node $N$ is labelled by $f_{\tau_T(N)}$.  In
fact, specifying $\tau_T$ corresponds to specifying amalgamation
maps in the free multicategory ${\mathbb F}$, and this defines the
amalgamation maps of $M_+$.

\end{enumerate}

Note that whereas in the symmetric case $\rho$ and $\tau$ may be
chosen freely for any given $T$, in this case precisely one
$\rho_T$ and $\tau_T$ is specified for each $T$. The source and
target of such an arrow $T$ are given by $s(T) = ( f_1, \ldots f_n
)$ and $t(T) = f \in a(M)$, the result of composing the $f_i$
according to their positions in $T$. Here, the tree $T$ may be
thought of as a combed tree as in the symmetric case, but with all
edges labelled by identities.

\begin{itemize}
\item Composition
\end{itemize}

When it can be defined, we have $T_1 \circ_m T_2 = T$ as follows:
\begin{enumerate}
\item $T$ is the combed labelled tree obtained from $(T_1,
\tau_{T_1})$ by replacing the node ${\tau_{T_1}}^{-1}(m)$ by the
combed tree $(T_2,\tau_{T_2})$, combing the tree and then
forgetting the twist at the top.

\item The amalgamation maps are defined to reorder the source
as necessary according to $\tau_{T_1}$, $\tau_{T_2}$ and $\tau_T$.

\end{enumerate}

\begin{itemize}
\item Identities: $1_f$ is the tree with one node, labelled by $f$.
\end{itemize}

This definition is easily seen to satisfy the axioms for a
generalised multicategory.  Note that a different choice of
amalgamation maps for ${\mathbb F}$ gives rise to different
bijections $\tau_T$ and hence different amalgamation maps in
$M_+$, resulting in an isomorphic slice multicategory.

\subsection{Comparison of slice}
\label{slicesymgen}

In this section we compare the slice constructions and make
precise the sense in which they correspond to one another. Recall
(section~\ref{tacxi}) that we have defined a functor
    \[\cat{GenMulticat} \map{\xi} \cat{TidySymMulticat}.\]
We now show that this functor `commutes' with slicing, up to
equivalence.

We will eventually prove (Corollary~\ref{mcatcord}) that for any
generalised multicategory $M$
    \[\xi(M_+) \simeq \xi(M)^+.\]
We prove this by constructing, for any  morphism of symmetric
multicategories $\phi:Q \longrightarrow \xi(M)$ a morphism
$\phi^+:Q^+ \longrightarrow \xi(M_+)$ such that
    \[\phi \mbox{ is an equivalence} \Rightarrow \phi^+ \mbox{ is an
    equivalence}.\]
The result then follows by considering the case $\phi = 1$.

We begin by constructing $\phi^+$. Recall
\[\begin{array}[t]{rcll}
    o(Q^+) & = & a(Q) & \\
    a(Q^+) & = & \{(T,\rho,\tau): &T \mbox { a labelled tree with
    $n$ nodes, $k$ leaves}\\
        &&& \rho \in {\mathbf S}_k,\\
        &&& \tau:\{\mbox{nodes of T}\} \map{\sim} [n]\\
        &&& \mbox{edges labelled by morphisms of }{\mathbb C}\}\\
    o(\xi(M_+)) & = & a(M) \\
    a(\xi(M_+)) & = & \{(T,\sigma)\ : &T \mbox{ a labelled tree
    with $n$ nodes} \\
        &&& \sigma\in {\mathbf S}_n \\
        &&& \mbox{edges labelled by identities}\}.
\end{array}\]

The idea is that given a way of composing arrows $f_1, \ldots,
f_n$ of $Q$ to an arrow $f$, we have a way of composing arrows
$g_1, \ldots, g_n$ of $M$ to an arrow $g$, where
    \begin{eqnarray*} \phi(f_i) & = & (g_i, \sigma_i) \\
    \mbox{and\ } \phi(f) & = & (g, \sigma).\end{eqnarray*}
Observe that since $\xi M$ is object-discrete, we have $\phi a =
1$ for all object-morphisms $a\in \bb{C}$.

So we define $\phi^+$ as follows:
\begin{itemize}
\item On objects: if $\phi(f) = (g,\sigma),\  g \in a(M)$ then put
$\phi^+(f) = g$.

\item On object-morphisms: since $\xi(M^+)$ is object-discrete,
we must have $\phi^+(\alpha) = 1$ for all object-morphisms
$\alpha$.

\item On arrows: put $\phi^+:(T,\rho,\tau)\longmapsto(\bar{T},\tau\circ
{\tau_{\bar{T}}}^{-1})$, where $\bar{T}$ is the labelled planar
tree obtained as follows.  Given a node with label $f$ say, and
$\phi(f) = (g,\sigma)$:
\begin{enumerate}
    \item replace the label with $g$
    \item `twist' the inputs of the node according to $\sigma$
    \item proceed similarly with all nodes, make all edge labels
    identities, then comb and ignore the twist at the top of the
    resulting tree (since the twist in $M_+$ is determined by the
    tree).
\end{enumerate}
\end{itemize}

For example, suppose $T$ is given by

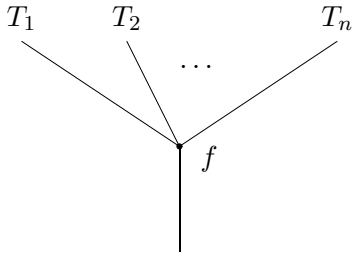
\setlength{\unitlength}{0.7mm}
\begin{center}
\begin{picture}(80,55)      

\put(40,0){\line(0,1){20}}   %
\put(40,20){\line(-3,2){30}} %
\put(40,20){\line(-1,2){10}} %
\put(40,20){\line(3,2){30}}  %
\put(40,20){\circle*{1}}     %

\put(10,43){\makebox[0pt]{$T_1$}}  
\put(30,43){\makebox[0pt]{$T_2$}}  
\put(70,43){\makebox[0pt]{$T_n$}}  

\put(44,20){\makebox(0,0)[tl]{$f$}} 

\put(40,35){$\ldots$}

\end{picture}
\end{center}where the $T_i$ are subtrees of $T$, and $\phi(f)=(g,\sigma)$.
Then steps (i) and (ii) above give

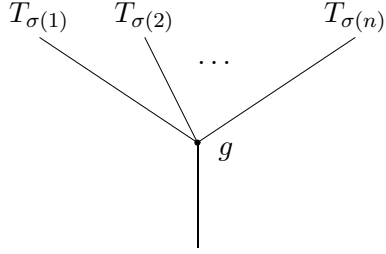
\setlength{\unitlength}{0.7mm}
\begin{center}
\begin{picture}(80,55)      

\put(40,0){\line(0,1){20}}   %
\put(40,20){\line(-3,2){30}} %
\put(40,20){\line(-1,2){10}} %
\put(40,20){\line(3,2){30}}  %
\put(40,20){\circle*{1}}     %

\put(10,43){\makebox[0pt]{$T_{\sigma(1)}$}}  
\put(30,43){\makebox[0pt]{$T_{\sigma(2)}$}}  
\put(70,43){\makebox[0pt]{$T_{\sigma(n)}$}}  

\put(44,20){\makebox(0,0)[tl]{$g$}} 

\put(40,35){$\ldots$}

\end{picture}
\end{center} and $\bar{T}$ is then defined inductively on the subtrees.  Node
$N$ in $\bar{T}$ is considered to be the image of node $N$ in $T$
under the operation $T \longrightarrow \bar{T}$.


Writing
\begin{eqnarray*}
    s(T,\rho,\tau) &=& (f_1,\ldots,f_n)\\
    \makebox(0,0)[br]{and \ \ }t(T,\rho,\tau) &=& f
\end{eqnarray*} we check that
\begin{eqnarray*}
    s(\phi^+(T,\rho,\tau)) &=& (\phi^+(f_1),\ldots,\phi^+(f_n))\\
     \makebox(0,0)[br]{and \ \ } t(\phi^+(T,\rho,\tau))
     &=&\phi^+(f).
\end{eqnarray*}
Writing $s(\bar{T},\tau\circ\tau_{\bar{T}}^{-1}) =
(g_1,\ldots,g_n)$ in $\xi(M)$, we have, in $M_+$
    \[s(\bar{T}) =
    (g_{\tau\circ\tau_{\bar{T}}^{-1}(1)},\ldots,
    g_{\tau\circ\tau_{\bar{T}}^{-1}(n)})\]
so node $N$ is labelled in $\bar{T}$ by
    \(g_{\tau\circ\tau_{\bar{T}}^{-1}(\tau_{\bar{T}}(N))} =
    g_{\tau(N)}\)
and in $T$ by
    \(f_{\tau(N)}.\)
So by definition of $\bar{T}$ we have
    \[\phi^+(f_{\tau(N)}) = g_{\tau(N)}\]
so $\phi^+(f_i) = g_i$ for each $i$ and
    \[s(\bar{T},\tau\circ\tau_{\bar{T}}^{-1}) =
    (\phi^+(f_1),\ldots,\phi^+(f_n))\]
as required.  Also, $t(\bar{T},\tau\circ\tau_{\bar{T}}^{-1}) =
\phi^+(f)$ by functoriality of $\phi$ and definition of
composition in $\xi(M)$.

We have shown that $\phi^+$ is functorial on the object-category
$o(Q^+)$; we need to check the remaining conditions for $\phi^+$
to be a morphism of symmetric multicategories. We may now assume
that all edge labels are identities since they all become
identities under the action of $\phi^+$.

\begin{itemize}
\item$\phi^+$ preserves identities:
\end{itemize}

$1_f \in a(Q^+)$ is $(T,\iota,\iota)$ where $T$ has one node,
labelled by $f$. So we have $\phi^+(1_f)=T$ where $T$ has one
node, labelled by $\phi^+(f)$, and $\phi^+(1_f)=1_{\phi^+}(f)$.

\begin{itemize}
\item $\phi^+$ preserves composition: We need to show
    \[\phi^+(\alpha\circ_m\beta)=\phi^+(\alpha)\circ_m\phi^+(\beta).\]
\end{itemize}

Now, the underlying trees are the same by functoriality of $\phi$,
the permutation of leaves is the same by coherence for
amalgamation maps of $M$, and the node ordering is the same by
definition of $\phi^+$.

\begin{itemize}
\item $\phi^+$ preserves symmetric action:
    \begin{eqnarray*}
    \phi^+((T,\rho,\tau)\sigma)) & = &
    \phi^+(T,\rho,\sigma^{-1}\tau)\\
        & = &
        (\bar{T},\sigma^{-1}\tau\circ{\tau_{\bar{T}}}^{-1})\\
        & = & (\bar{T},\tau\circ{\tau_{\bar{T}}}^{-1})\sigma\\
        & = & (\phi^+(T,\rho,\tau))\sigma.\\
    \end{eqnarray*}
\end{itemize}
So $\phi^+$ is a morphism of symmetric multicategories.

\begin{proposition} \label{mcatpropb} Let $Q$ be a symmetric
multicategory, $M$ a generalised multicategory and
$\phi:Q\longrightarrow\xi(M)$ a morphism of symmetric
multicategories. If $\phi$ is an equivalence then $\phi^+$ is an
equivalence.
\end{proposition}


This enables us to prove the following proposition:

\begin{proposition} \label{mcatpropf} If $Q$ is tidy then $Q^+$ is tidy.
\end{proposition}

\noindent {\bfseries Proof of Proposition \ref{mcatpropb}.  }
First we observe that given any such morphism $\phi$, $Q$ is
freely symmetric:
    \begin{eqnarray*}
    \alpha\sigma = \alpha & \Rightarrow & \phi(\alpha\sigma) =
    \phi(\alpha)\sigma = \phi(\alpha) \ \in \xi(M) \\
    & \Rightarrow & \sigma = \iota,
    \end{eqnarray*}
the second implication following from $\xi(M)$ being freely
symmetric.

Now, given that $\phi$ is full, faithful and essentially
surjective on the category of objects, and full and faithful, we
prove the proposition in the following steps:

\renewcommand{\labelenumi}{\roman{enumi})}
\begin{enumerate}
\item $\phi^+$ is surjective on objects
\item $\phi^+$ is full on the category of objects
\item $\phi^+$ is faithful on the category of objects
\item $\phi^+$ is full
\item $\phi^+$ is faithful
\end{enumerate}

\begin{prfof}{(i)}  Recall the action of $\phi^+$ on objects: let $f\in o(Q^+)
= a(Q)$ with $\phi(f)=(g,\sigma)$ then $\phi^+:f\longmapsto g$.
Now, given any $g\in o(\xi(M_+)) = a(M)$, we have $(g,\iota)\in
a(\xi( M))$. $\phi$ is full and surjective, so there exists $f\in
a(Q)$ such that $\phi(f)=(g,\sigma)$ and $\phi^+(f)=g$.
\end{prfof}

\begin{prfof}{(ii)}$\xi(M_+)$ is object-discrete so we only need to show that
if $\phi^+(f_1)=\phi^+(f_2)$ then there is a morphism $f_1
\longrightarrow f_2$ in $o(Q^+)$. Now
    \begin{eqnarray*}
    \phi^+(f_1)=\phi^+(f_2) \ \Rightarrow \ \phi(f_1) & = &
    \phi(f_2)\sigma \mbox{ \ for some permutation } \sigma\\
        & = & \phi(f_2 \sigma).
    \end{eqnarray*}
Suppose
\[\begin{array}{rcccc}
    f_1 & : & a_1, \ldots, a_n & \lra & a
    \\
   \mbox{and \ } f_2 \sigma & : & b_1, \ldots, b_n & \lra & b. \\
    \end{array}\]
Then we must have $\phi(a_i)=\phi(b_i)$ for all $i$, and
$\phi(a)=\phi(b)$. So there exist morphisms
\[\begin{array}{rcccc}
    g_i & : & b_i & \lra & a_i
    \\
   \mbox{and \ } g & : & a & \lra & b \\
    \end{array}\]
and we have
    \[f_2\sigma = g \circ f_1 \circ
    (g_1,\ldots,g_n)\]
giving a morphism $f_1 \lra f_2$ as required. \end{prfof}

\begin{prfof}{(iii)} An arrow $\alpha:f_1 \longrightarrow f_2$ is uniquely of the
form $(\sigma,g_1,\ldots,g_n;g)$ with
\[\begin{array}{rcccc}
    g_i & : & s(f_2)_{\sigma(i)} & \lra & s(f_1)_i
    \\
   \mbox{and \ } g & : & t(f_1) & \lra & t(f_2) \\
    \end{array}\]
as arrows of ${\mathbb C}$. Since $\phi$ is faithful on the
category of objects and $\xi(M)$ is object-discrete, there can
only be one such map. \end{prfof}


\begin{prfof}{(iv)} Given $f_1,\ldots,f_n,f \in o(Q^+)$ and
    \[(T,\sigma):(\phi^+(f_1),\ldots,\phi^+(f_n))
    \longrightarrow\phi^+(f) \in \xi(M_+)\]
we seek
    \[(T',\rho,\tau):(f_1,\ldots,f_n) \longrightarrow f \in
    Q^+\]
such that
    \[\phi^+(T',\rho,\tau)=(T,\sigma)\]
i.e. such that $\bar{T'}=T$ and
$\tau\circ{\tau_{\bar{T}}}^{-1}=\sigma$.

Write $\phi(f)=(g,\alpha)$ and for each $i$, $\phi(f_i) =
(g_i,\alpha_i)$.  Then $\phi^+(f_i)=g_i$ and $\phi^+(f)=g$.
$(T,\sigma)$ is a configuration for composing the $g_i$ to yield
$g$, so we certainly have a configuration for composing the
$(g_i,\alpha_i)$ to yield $g_i$ as follows: replace node label
$g_i$ by $(g_i,\alpha_i)$ and insert a twist ${\alpha_i}^{-1}$
above the node, then comb and add the necessary twist at the top.

This gives a configuration for composing the $f_i$ as follows.  We
have
    \[t(g_i,\alpha_i) = s(g_k,\alpha_k)_m \Rightarrow \phi(t(f_i))
    = \phi(s(f_k)_m).\]
Now $\phi$ is faithful on the category of objects, so there exists
a morphism
    \[t(f_i) \longrightarrow s(f_k)_m\]
and we label the edge joining $t(f_i)$ and $s(f_i)_m$ with this
object-morphism.  So this gives a configuration for composing the
$f_i$, to yield $h$, say, with $\phi(h) = \phi(f)$. That is, we
have a morphism
    \[(f_1,\ldots,f_n) \stackrel{\theta}{\longrightarrow} h\]
such that $\phi^+(\theta) = (T,\sigma)$.

Now $\phi$ is full on the category of objects, so if $\phi(h) =
\phi(f)$ then there is a morphism $\alpha:h \longrightarrow f$ in
$o(Q^+)$.  So we have
    \[(f_1,\ldots,f_n)\stackrel{\theta}{\longrightarrow} h
    \stackrel{\iota(\alpha)}{\longrightarrow}f\]
and $\phi^+(\iota(\alpha))$ is the identity since $\xi(M_+)$ is
object-discrete.  So
    \[\phi^+(\iota(\alpha)\circ\theta) = \phi^+(\theta) =
    (T,\sigma)\]
as required.  \end{prfof}

\begin{prfof}{(v)} Suppose $\phi^+(\alpha)=\phi^+(\beta)$.  Then, writing
\[\begin{array}{rcccc}
    \alpha = (T_1,\rho_1,\tau_1) & : & (f_1,\ldots,f_n) & \lra & f
    \\
    \beta = (T_2,\rho_2,\tau_2) & : & (f_1,\ldots,f_n) & \lra  & f
\end{array}\]
we have $\bar{T_1}=\bar{T_2}=\bar{T}$, say, and
$\tau_1\circ{\tau_{\bar{T_1}}}^{-1} =
\tau_2\circ{\tau_{\bar{T_2}}}^{-1}$ so $\tau_1=\tau_2$. So given
any node $N$ in $\bar{T}$, its pre-image in $T_1$ has the same
label $f_i$ as its pre-image in $T_2$.  The same is true of edge
labels, since $\phi$ is faithful on the category of objects.

Then the tree $T_1$ may be obtained from $\bar{T}$ as follows.
Suppose $\phi(f_i) = (g_i,\sigma)$ and $\phi(f)=g$.  Then for the
node labelled by $g_i$, apply the twist $\sigma^{-1}$ to the edges
above it, and then relabel the node with $f_i$. This process may
also be applied to obtain the tree $T_2$.  Since the process is
the same in both cases, we have $T_1=T_2=T$, say.

Finally, suppose $f'$ is the arrow obtained from composing
according to $T$.  Then by the action of $\alpha$, $f=f'\rho_1$,
and by the action of $\beta$, $f=f'\rho_2$. Then, since $Q$ is
freely symmetric, $\rho_1=\rho_2$, so $\alpha=\beta$ as required.
\end{prfof}

\begin{prfof}{Proposition \ref{mcatpropf}}  Given a tidy symmetric
multicategory $Q$ we need to show that $Q^+$ is also tidy.

Recall (Lemma~\ref{tidylemma}) that a \sm\ $Q$ is tidy if and only
if it is equivalent to one in the image of $\xi$, $\xi M$ say,
with equivalence given by
    \[\phi: Q \lra \xi (M).\]
Then by Proposition~\ref{mcatpropb} $\phi^+$ is an equivalence
    \[\phi^+: Q^+ \lra \xi(M_+)\]
so $Q^+$ is tidy as required.
\end{prfof}

\begin{corollary} \label{mcatcord} Let $M$ be a generalised
multicategory.  Then \[\xi(M)^+ \simeq \xi(M_+)\] as symmetric
multicategories with a category of objects.
\end{corollary}

\begin{prf}  Put $Q=\xi(M)$, $\phi=1$ in Proposition~\ref{mcatpropb}.
\end{prf}

We are now ready to give the construction of opetopes.  
\label{tacopetopes}

\subsection{Opetopes}
\label{tacopes}

For any symmetric multicategory $Q$ we write
    \[Q^{k+} = \left\{\begin{array} {l@{\extracolsep{2em}}l}
    Q & k=0 \\
    {(Q^{(k-1)+})}^+ & k \ge 1\end{array} \right. \]
Let $I$ be the symmetric multicategory with precisely one object,
precisely one (identity) object-morphism, and precisely one
(identity) arrow.  A {\em $k$-dimensional opetope}, or simply {\em
$k$-opetope}, is defined in \cite{bd1} to be an object of
$I^{k+}$. We write $\bb{C}_k = o(I^{k+})$, the category of
$k$-opetopes.

\subsection{Multitopes}
\label{mtopes}

Multitopes are defined in \cite{hmp1} using the multicategory of
function replacement.  We give the same construction here, but
state it in the language of slicing; this makes the analogy with
Section~\ref{tacopes} clear.

For any generalised multicategory $M$ we write
    \[M_{k+} = \left\{\begin{array} {l@{\extracolsep{2em}}l}
    M & k=0 \\
    {(M_{(k-1)+})}_+ & k \ge 1 \end{array} \right.\]
Let $J$ be the generalised multicategory with precisely one object
and precisely one (identity) morphism. Then a {\em$k$-multitope}
is defined to be an object of $J_{k+}$.  We write $P_k =
o(J_{k+})$, the set of $k$-multitopes; we will also regard this as
a discrete category.

\subsection{Comparison of opetopes and multitopes} \label{opecomp}

In this section we compare the construction of opetopes and
multitopes, applying the results we have already established.

\begin{proposition} For each $k \ge 0$
    \[\xi(J_{k+}) \simeq I^{k+}.\]
\end{proposition}

\begin{prf}  By induction.  First observe that $\xi(J) \cong I$ and write $\phi$ for this
isomorphism. So for each $k \ge 0$ we have
    \[\phi^{k+}:I^{k+} \longrightarrow \xi(J_{k+}),\]
where
     \[\phi^{k+} = \left\{\begin{array} {l@{\extracolsep{1.5cm}}l}
    \phi & k=0 \\
    (\phi^{(k-1)+})^+ & k \ge 1
\end{array} \right. \]

Now $I$ is (trivially) tidy, so by Proposition~\ref{mcatpropf},
$I^{k+}$ is tidy for each $k \ge 0$. So by
Proposition~\ref{mcatpropb}, $\phi^{k+}$ is an equivalence for all
$k \ge 0$. \end{prf}

Then on objects, the above equivalence gives the following result.

\begin{corollary} For each $k \geq 0$
    \[P_k \simeq {\mathbb C}_k.\]
\end{corollary}

This results shows that `opetopes and multitopes are the same up to
isomorphism'.

\addcontentsline{toc}{section}{References}
\bibliography{bib0209}

\nocite{bd2}

\nocite{hmp2} \nocite{hmp3} \nocite{hmp4}

\nocite{bae1}

\nocite{ks1}

\nocite{sim1}

\nocite{ben1}

\end{document}